\documentclass[10pt,reqno]{amsart}
\usepackage{chemarrow}
\usepackage{amsmath}
\usepackage{amssymb}
\usepackage{amsthm}
\usepackage{color}
\usepackage{graphicx}
\usepackage{psfrag}
\usepackage{marginnote}
\usepackage{mathrsfs}
\usepackage{enumerate}
\usepackage[babel,german=quotes]{csquotes}

\renewcommand{\ge}{\geqslant}
\renewcommand{\geq}{\geqslant}

\renewcommand{\le}{\leqslant}
\renewcommand{\leq}{\leqslant}

\newcommand{\Nb}{{\mathbb N}}
\newcommand{\Rb}{{\mathbb R}}

\newcommand{\sign}{\operatorname{sign}} 

\makeatletter
\renewcommand*\env@matrix[1][\arraystretch]{%
  \edef\arraystretch{#1}%
  \hskip -\arraycolsep
  \let\@ifnextchar\new@ifnextchar
  \array{*\c@MaxMatrixCols c}}
\makeatother

\newenvironment{Proof}{\removelastskip \vskip12pt plus 1pt \noindent
{\em Proof\/}\rm }{\hfill$\square$ \vskip12pt plus 1pt}

\newtheorem{theorem}{Theorem}[section]
\newtheorem{lemma}[theorem]{Lemma}
\newtheorem{proposition}[theorem]{Proposition}
\newtheorem{definition}[theorem]{Definition}

\theoremstyle{definition}
\theoremstyle{remark}

\usepackage{showlabels}
\usepackage{lscape}
\usepackage{hyperref}

\begin{document}
\title[Generalized Exchange Driven Model]{The continuous version of the generalized exchange-driven growth model}
\date{November 30, 2023}

\author{P. K. Barik}
\address[P. K. Barik]{Instituto de Matemáticas,
Universidad de Granada, Rector López Argüet, 
S/N, 18001, Granada, Spain, 
and
Departamento de Matemática Aplicada,
Universidad de Granada,
Avenida de Fuentenueva S/N, 18071, Granada, Spain.}  \email{barik@ugr.es}

\author{F. P. da Costa}
\address[F. P. da Costa]{Univ. Aberta, Dep. of Sciences and Technology,
  Rua da Escola Polit\'ecnica 141-7, P-1269-001 Lisboa, Portugal, and
  Univ. Lisboa, Instituto Superior T\'ecnico, Centre for Mathematical
  Analysis, Geometry and Dynamical Systems, Av. Rovisco Pais,
  P-1049-001 Lisboa, Portugal.}  \email{fcosta@uab.pt}

 \author{J. T. Pinto} 
\address[J. T. Pinto]{Univ. Lisboa, Instituto Superior T\'ecnico, Dep. of Mathematics
and Centre for Mathematical
Analysis, Geometry and Dynamical Systems, Av. Rovisco Pais,
P-1049-001 Lisboa, Portugal.}  \email{jpinto@tecnico.ulisboa.pt}

\author{R. Sasportes} 
\address[R. Sasportes]{Univ. Aberta, Dep. of Sciences and Technology,
Rua da Escola Polit\'ecnica 141-7, P-1269-001 Lisboa, Portugal, and
Univ. Lisboa, Instituto Superior T\'ecnico, Centre for Mathematical
Analysis, Geometry and Dynamical Systems, Av. Rovisco Pais,
P-1049-001 Lisboa, Portugal.}  \email{rafael.sasportes@uab.pt}

\begin{abstract}
In this article, we discuss the continuous version of the generalized exchange-driven growth model which is a variant of the coagulation model in which a smaller size particle is detached from a bigger one and merges with another particle. This new model is a continuous extension of the generalized exchange-driven growth model originally formulated in a discrete context \cite{daCosta:2023}. In this work, we examine the existence of weak solutions to the continuous version of the generalized exchange-driven growth model under a suitable reaction rate. Under an additional condition on the reaction rates, a uniqueness result is established. Finally, we prove that solutions satisfy the mass-conserving property and the conservation of the total number of particles for coagulation rates with linear bounds.
\end{abstract}

\maketitle
\section{Introduction} \label{sec:intro}
Over the past five decades, coagulation-fragmentation type models have become the focus of a substantial
the volume of research literature in the mathematical community (e.g.: see \cite{bll,dc15} and references therein).
They have applications in the study of a large range of physical, biological, ecological, economical, and social phenomena, such as, for example:
formation of stars and planets \cite{safronov}, deposition processes in solid state physics \cite{mulheran}, 
aerosol dynamics \cite{marina}, coagulation of red blood cells \cite{guy}, grouping behavior of 
animals \cite{degond}, merger of companies \cite{banakar}, or in some toy models of social behavior \cite{ispolatov}.

\medskip

These models describe the time evolution of a system of entities, usually called \emph{particles} or \emph{clusters}, 
which are characterized by a parameter(s) describing some relevant aspect, such as their sizes or masses. Particles can 
grow (i.e., the parameter can increase) or shrink (the parameter decreases) due to the action of opposing mechanisms, 
typically called \emph{coagulation} (the growth process) and \emph{fragmentation} (the
shrinking mechanism). 
To elaborate further, coagulation processes occur when two particles interact and stick together producing a larger particle. 
Conversely, in the fragmentation process, a bigger particle breaks into two or more than two daughter particles, due
to external factors or to spontaneous internal ones.

\medskip

Coagulation-fragmentation models describe these populations of particles and their time evolution by means of an infinite system of  ordinary or partial differential equations  or by partial integro-differential equations, depending on the application context of the given problem. If particle sizes are discrete they are usually modelled by non-negative integers, and the models consist of an infinite system of differential equations.  Conversely, if particle sizes can vary continuously, they are usually described by non-negative real numbers and the  models consist of a partial integro-differential equation. 

\medskip

\begin{figure}[h]
	\includegraphics[scale=0.9]{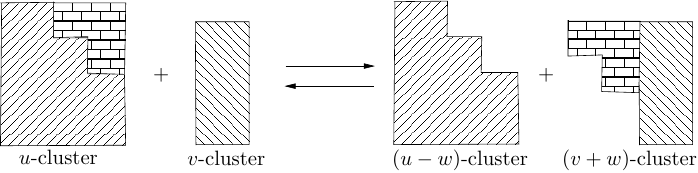}\label{fig1}
	\caption{Scheme of the generalized exchange-drive mechanism. In the forward reaction scheme a $u$-cluster 
	 comes into contact 
with a $v$-cluster, 
 and a part of the $u$-cluster with mass $w$ 
 becomes attached to the $v$-cluster. Analogously for the backward reaction.}
\end{figure}

Several variants of these basic processes can be relevant resulting in a plethora of different specific models. In this paper we  discuss the continuous version of the generalized exchange-driven growth model which is a generalization of the exchange-driven system \cite{esen1,esen2,esenvel} to both continuous cluster sizes and to the possibility of exchanging chunks of any mass between two clusters. To clarify the phenomena to be studied, let us consider a system of particle clusters with masses $u\in\Rb_{>0}$ and denote a cluster of mass $u$, or $u$-cluster, by $\langle u \rangle$. Now consider that when two clusters meet a chunk of one of the clusters becomes attached to the other, more precisely let us imagine a situation where two particles, $\langle u \rangle$ and $\langle v \rangle$, come into contact, and this interaction triggers the detachment of a smaller part of $\langle u \rangle$, with mass $w \le u $, from cluster $\langle u \rangle$, immediately followed by its attachment to the cluster $\langle v \rangle$, by a reaction taking place a rate denoted by $\widetilde{\mathcal{A}}(u, v; w)$. This is schematically represented by
\begin{equation}
\langle u \rangle + \langle v\rangle \xrightarrow[{}]{\widetilde{\mathcal{A}}(u, v; w)}  \langle u-w\rangle + \langle v+w\rangle.  \label{kinetic}
\end{equation}
and is illustrated in Figure~\ref{fig1}.

%
%
%
%
%
%

\medskip

Clearly, here we assume $u, v, w\in \mathbb{R}$ with $u > 0$, $0 < w \leq u$, and $v\geq 0$. This is because the case of $w=0$ corresponds to a lack of reaction, while $w>u$ lacks physical meaning, as it is impossible to detach a portion larger than the whole particle size from it. However, we will allow $w=u$, implying the entire particle of size $u$ becomes attached to the particle of size $v$ in the reaction:
\[
\langle u \rangle + \langle v\rangle \xrightarrow[{}]{\widetilde{\mathcal{A}}(u, v; u)}  \langle 0\rangle + \langle v+u\rangle.
\]
This closely resembles the conventional coagulation reaction, except now we introduce the concept of a ``void'' cluster $\langle 0\rangle$.
In a similar fashion, if we consider $v=0$ then the scheme transforms into:
\[
\langle u\rangle + \langle 0\rangle \xrightarrow[{}]{\widetilde{\mathcal{A}}(u, 0; w)} \langle u-w\rangle + \langle w\rangle.
\]
This is a type of fragmentation of the non-void particle size. Despite the introduction of a nonlinear contribution due to the presence of the void cluster in the law of mass action, when we envision the system within an infinite bath of void particles at a constant concentration, this reaction precisely corresponds to the binary fragmentation reaction explored in coagulation-fragmentation studies.

\medskip

Under reasonably diluted conditions (i.e., at low cluster concentrations) one can assume the mass action law to be valid. If the equations are chemically elementary  the process schematically represented by \eqref{kinetic} and by Figure~\ref{fig1} can be modelled by a system of chemical 
kinetics differential equations describing the time evolution of each type of cluster in which the rate of the chemical equation \eqref{kinetic} is given by
\[
\widetilde{\mathcal{A}}(u, v; w) := \mathcal{A}(u, v; w)\zeta(t, u)\zeta(t, v),
\]
where $\zeta(t, u)$ is the concentration of $u$-clusters at time $t$, and $ \mathcal{A}(u, v; w)$ is a parameter independent of the concentrations of the clusters involved in the reaction, called the rate constant, or rate coefficient, of the reaction. Next, we clarify how the full system of differential equations modelling the evolution of the different clusters comes about.

\medskip

With the above assumptions the rate of change of the concentration of $x$-clusters, $\zeta(t, x),$ is determined by several terms, that can essentially be grouped into two classes that we now detail: Let us start by observing that $\langle x \rangle$ can be produced from reactions where a bigger cluster sheds an appropriately sized part and, in a corresponding reverse reaction, it can be destroyed by becoming bigger due to the capture of part of another cluster; schematically we have:
\begin{equation}
\langle x+z \rangle + \langle y \rangle \;\autorightleftharpoons{$\widetilde{\mathcal{A}}(x+z, y; z)$}{$\widetilde{\mathcal{A}}(y+z, x; z)$}\;  \langle x\rangle + \langle y+z\rangle.
\label{reaction1}
\end{equation}
Clearly, in principle, both $y$ and $z$ can be any positive real. So, by the mass action law, 
the contributions of these reactions for the rate of change of $\zeta(t,x)$ result in
\begin{equation}
 \mathcal{B}_1( \zeta(t, \cdot) )(x) :=  \int_0^{\infty}  \int_0^{\infty} \mathcal{A}(x+z, y; z)  \ \zeta(t, x+z) \ \zeta(t, y) \ dy  \ dz;		\label{B1}
\end{equation}
for its increase due to the forward reaction in \eqref{reaction1}, whereas the destruction of clusters $\langle x\rangle$ through the opposite reaction
contribute to the term
\begin{equation}
 \mathcal{D}_1( \zeta(t, \cdot) )(x) :=  - \int_0^{\infty}  \int_0^{\infty} \mathcal{A}(y+z, x; z) \  \zeta(t, y+z) \ \zeta(t, x) \ dy \ dz.		\label{D1a}
\end{equation}
It will be convenient for the computations to change variables $y\mapsto \tilde{y}:=y+z$ in \eqref{D1a} and write it as 
\begin{equation}
 \mathcal{D}_1( \zeta(t, \cdot) )(x) :=  - \int_0^{\infty}  \int_z^{\infty} \mathcal{A}(y, x; z) \  \zeta(t, y) \ \zeta(t, x) \ dy \ dz,		\label{D1}
\end{equation}
where we dropped the tilde.

\medskip

Now, $\langle x\rangle$ can also be produced by reactions involving the capture, by a smaller cluster, of an appropriately sized part of another cluster, as in the forward reaction in the scheme
\begin{equation}
\langle x-z \rangle + \langle y \rangle \;\autorightleftharpoons{$\widetilde{\mathcal{A}}(y, x-z; z)$}{$\widetilde{\mathcal{A}}(x, y-z; z)$}\;  \langle x\rangle + \langle y-z\rangle.
\label{reaction2}
\end{equation}
which, like in \eqref{reaction1}, also contemplates the reverse reaction. The contributions of these reactions to the overall rate of change of $\zeta(t,x)$ are, again due to the mass action law,
\begin{equation}
 \mathcal{B}_2( \zeta(t, \cdot) )(x) :=  \int_0^{x} \int_z^{\infty}   \mathcal{A}(y, x-z; z)  \ \zeta(t, y) \ \zeta(t, x-z) \ dy \ dz,	\label{B2}
\end{equation}
for the forward reaction in \eqref{reaction2}, and
\begin{equation}
 \mathcal{D}_2( \zeta(t, \cdot) )(x) := - \int_0^{x} \int_z^{\infty}   \mathcal{A}(x, y-z; z) \  \zeta(t, x) \ \zeta(t, y-z) \ dy \ dz		\label{D2a}
\end{equation}
for the backward reaction. Again as in \eqref{D1a} above, it will be convenient to change variables $y\mapsto \tilde{y}:=y-z$ in
\eqref{D2a} and write $\mathcal{D}_2( \zeta(t, \cdot) )(x)$ as
\begin{equation}
 \mathcal{D}_2( \zeta(t, \cdot) )(x) := - \int_0^{x} \int_0^{\infty}   \mathcal{A}(x, y; z) \  \zeta(t, x) \ \zeta(t, y) \ dy \ dz		\label{D2}
\end{equation}
where we, again, dropped the tilde.

\medskip

\noindent  Thus, the overall rate of change of $\zeta(t, x)$, the concentration of particles of mass $x \in \mathbb{R}_{>0}:= (0, \infty)$ at time $t \ge 0$, can be expressed by the following  \emph{continuous version of the generalized exchange-driven growth system\/} (CGEDG)
\begin{align}\label{GEDM}
\frac{\partial \zeta(t, x) }{\partial t}  =  \Bigl(\mathcal{B}_1 ( \zeta(t, \cdot) )(x)  + \mathcal{D}_1 ( \zeta(t, \cdot) )(x)\Bigr)  
+ \Bigl(\mathcal{B}_2 ( \zeta(t, \cdot) )(x) + \mathcal{D}_2 ( \zeta(t, \cdot) )(x)\Bigr).
\end{align}
We will consider initial value problems for \eqref{GEDM}, with initial conditions
\begin{align}\label{Initialdata}
  \zeta(0, x) = \zeta^{in}(x) \ge 0\ a.e.
\end{align}

\medskip

Next, we proceed to define the moments of particles for the generalized exchange-driven growth equation \eqref{GEDM}. For $r \ge 0$, $\mathcal{M}_r(t)$ represents the $r^{th}$  moment of particles of size $x$ at time $t$ and if it exists, it can be defined as follows:
 \begin{align}\label{Totalmass}
\mathcal{M}_r(t) =\mathcal{M}_r(\zeta(t)):=\int_0^{\infty} x^r \ \zeta(t, x)\ dx,  \ \ t \ge 0.
\end{align}
The zeroth and first moments correspond to specific physical quantities. When $r=0$, $\mathcal{M}_0(t)$ corresponds to the total number of particles, and for $r=1$, $\mathcal{M}_1(t)$ represents the total mass of particles within the system.

\medskip

The primary goal of this article is to discuss the fundamental properties of the model \eqref{GEDM}--\eqref{Initialdata}, including its 
existence, uniqueness, and physical characteristics such as mass conservation and the preservation of the total number of particles. 
Before delving into these topics, we will begin by providing a concise overview of existing articles related to this class of models. 
The mathematically rigorous study of the discrete exchange-driven growth model 
that considers only clusters with discrete masses $i\in\{0, 1, \ldots\}$ and allows 
only for the exchange of a particle of mass $1$ when two clusters come into interaction, was
the subject of \cite{esen1, esenvel, schl1}, where 
existence, uniqueness, and physical properties such as mass conservation of solutions were proved.
In these papers, some aspects of the long-time behaviour of solutions were also investigated. 
Moreover, in \cite{Eichenberg:2021}  there is a detailed discussion of 
self-similar solutions. Recently this model was generalized for the 
possibility of exchanging particles of higher masses: in \cite{daCosta:2023} we started the 
study of existence, uniqueness, and physical
properties of the discrete version of this generalized discrete exchange-driven growth model which is the discrete version 
of \eqref{GEDM}--\eqref{Initialdata}.

\medskip

The content of the present paper is the following: in Section 2 we introduce essential notations, the key theorems to be proved below, and some preliminary results. In Section 3 we establish the existence of weak solutions for both sublinear and nearly quadratic coagulation rates. Section 4 is dedicated to the proof of the uniqueness of solutions under specific conditions on the coagulation coefficients. In the final section, we prove that weak solutions preserve mass and conserve the total number of particles.

 \section{Statement of the problem and results}\label{sec:results}
 Here we denote $L^1_{0, 1}(\mathbb{R}_{>0}):= L^1(\mathbb{R}_{>0}, (1+x)dx )$ the space of integrable function with respect to measure $(1+x)dx$. Now, we denote the set of all continuous functions from $[0, \infty)$ into $L^1_{0, 1}(\mathbb{R}_{>0})$ equipped with its weak topology by $\mathcal{C}( [0, \infty),    w\text{-}L^1_{0, 1}(\mathbb{R}_{>0}))$ (set of all weakly continuous functions). More precisely, a function 
$g \in \mathcal{C}([0, \infty), w\text{-}L^1_{0, 1}(\mathbb{R}_{>0}))$ if 
\begin{align*}
t\longmapsto \int_0^{\infty}  (1+x) \psi(x) g(t, x) dx 
\end{align*}
is a continuous function on  $[0, \infty )$  for all   $\psi \in L^{\infty}( \mathbb{R}_{>0} )$. When no danger of
misunderstanding arises, to simplify notation we write $L_{0,1}^1$ instead of $L^1_{0, 1}(\mathbb{R}_{>0})$ and denote
simply by $\|\cdot\|$ the standard norm in  $L_{0,1}^1$, instead of $\|\cdot\|_{L_{0,1}^1},$ and analogously for the $L^\infty$
space and standard norm.

Let us now state the definition of the solution to \eqref{GEDM}--\eqref{Initialdata} in the following sense:
\begin{definition}\label{Definition of Soln}
Let $0 \le \zeta^{in} \in L^1_{0, 1}$. We say that $\zeta $ is a weak solution to \eqref{GEDM}--\eqref{Initialdata}  if it fulfills the conditions:\\

$(a)$ $0 \le \zeta  \in \mathcal{C}( [0, \infty),    w\text{-}L^1(\mathbb{R}_{>0})) \cap L^{\infty} (0, T; L^1_{0, 1})$,\\

$(b)$ 
\begin{align}\label{Weak form}
 \int_0^{\infty} \!\!\omega(x) \bigl[ \zeta(t, x) - \zeta^{in}(x) \bigr] dx 
 = \int_0^t \int_0^{\infty}\! \int_0^{\infty}\! \int_0^{x}  \tilde{\omega}(x, y, z) \ \mathcal{B}(x, y; z, \zeta(s) )   \ dz  dy dx ds,
\end{align} 
for any smooth test function $\omega \in L^{\infty}(\mathbb{R}_{>0})$, where
\begin{align*}
  \tilde{\omega}(x, y, z) :=   \omega(y+z) + \omega(x-z)  - \omega(x) - \omega(y),
\end{align*}
and 
\begin{align*}
\mathcal{B}(x, y; z, \zeta(t)) :=    \mathcal{A}(x, y; z)  \ \zeta(t, x) \ \zeta(t, y).
\end{align*}
\end{definition}

We assume that the coagulation kernel $\mathcal{A}$ is a non-negative measurable function
satisfying the following \emph{physical assumption}: 
\begin{equation}\label{physical}
\mathcal{A}(x, y; z) = 0 \;\,\text{if}\;\, z > x.
\end{equation}
In addition, 
it also satisfies the following class of growth conditions for some non-negative constant $A$:\\
Class-I: Sum kernel:
 \begin{equation}\label{SumK}
\mathcal{A}(x, y; z) \le \begin{cases}
A \ \varphi(z),\       & \text{if}\ x+y < 1,\ \\
A\ (x+y)  \varphi(z),\ &  \text{if}\ x+y > 1,
\end{cases}
\end{equation}
where $\varphi(z)$ is an integrable function such that
\begin{equation}\label{PhiIntegrable}
\varphi(z) \in L^1_{0, 1}.
\end{equation} 
or\\
Class-II: Product kernel: for some function $\eta: (0, \infty) \rightarrow [1, \infty)$
 \begin{equation}\label{ProductK}
\mathcal{A}(x, y; z) \le \begin{cases}
A \ \varphi(z),\ & \text{if}\ (x, y) \in (0, 1)^2,\ \\
A \ \eta(y)  \varphi(z), \ &  \text{if}\ (x, y) \in (0, 1) \times (1, \infty), \  \\
A \ \eta(x)  \varphi(z), \ &  \text{if}\ (x, y) \in  (1, \infty) \times (0, 1),\ \\
A \ \eta(x) \eta(y)  \varphi(z), \  &  \text{if}\ (x, y) \in (1, \infty)^2.
\end{cases}
\end{equation}
where $\eta$ satisfies the following conditions:
\begin{align}\label{RateCoag}
\eta^{\ast} :=\sup_{x \ge 1} \frac{\eta(x)}{(1+x)} < \infty\ \ \text{and} \ \lim_{x \to \infty} \frac{\eta(x)}{x} =0.  
\end{align}
From \eqref{SumK}, \eqref{ProductK} and \eqref{RateCoag}, we infer that
\begin{align}\label{EstimateCoagRate}
\mathcal{A}(x, y; z) \le A (1+x)(1+y) \varphi(z).    
\end{align}
Further, we assume an additional assumption on the kernel $\mathcal{A}$, i.e.,
\begin{align}\label{Sobolevbound}
   \mathcal{A} \in W^{1, \infty}_{loc}( \mathbb{R}_{>0} \times \mathbb{R}_{>0} \times \mathbb{R}_{>0}  ),
\end{align}
and also
\begin{align}\label{boundderivative}
\frac{ \partial} {\partial x} \mathcal{A}(x, y; z) \le A(1+y) \varphi(z)\ \ \text{and}\  \frac{ \partial} {\partial y} \mathcal{A}(x, y; z) \le A(1+x) \varphi(z).
\end{align}
Finally, we assume that the initial data $\zeta^{in}$ satisfies the following condition:
\begin{equation}\label{ConditionIntial}
\zeta^{in} \in L^1_{0, 1}.
\end{equation}
The cornerstone of our analysis in this article relies on a monotonically increasing convex function and its derivative, which is a concave function. 
We shall initially consider two non-negative, convex, and non-decreasing functions $\sigma_1$ and $\sigma_2 \in \mathcal{C}^{2}([0, \infty))$ 
such that
 \begin{align}\label{Convexp0}
    \sigma_j(0)=\sigma_j'(0)=0\ \text{and}\ \sigma_j' \ \text{is concave};
 \end{align}
  \begin{align}\label{Convexp1}
   \lim_{x \to \infty} \sigma_j'(x) =\lim_{x \to \infty} \frac{ \sigma_j(x)}{x}=\infty;
 \end{align}
 for $j=1, 2$.
 Finally, we assume additional integrability conditions on $\zeta^{in}$ and $\varphi$ by applying a refined version of the de la Vall\'{e}e Poussin theorem see  \cite[Theorem~2.8]{Laurencot:2015} which are described below:
 \begin{align}\label{Convexp2}
\Gamma_1 := \int_0^{\infty}\sigma_1(x) \zeta^{in}(x) dx < \infty,~~\text{and}~~\Gamma_2 :=\int_0^{\infty}{\sigma_2(\zeta^{in}(x))} dx <\infty,
\end{align}
and
\begin{align}\label{Convexp3}
\Gamma_3 := \int_0^{\infty} \sigma_1(x) \varphi(x) dx <\infty,~~\text{and}~~\Gamma_4 :=\int_0^{\infty}{\sigma_2( \varphi(x))} dx < \infty.
\end{align}
Let us revisit some supplementary properties of the convex functions $\sigma_1$ and $\sigma_2$ as outlined in \cite{Laurencot:2015}: 
\begin{lemma} Consider $\sigma_1$, $\sigma_2$ in $\mathcal{C}^2{ (\mathbb{R}_{>0})  }$ such that $\sigma'_1$ and $\sigma'_2$ are concave. Then, the following inequalities hold:
\begin{equation}\label{Convexp4}
\hspace{-5cm} \sigma_2(x)\le x \sigma'_2(x)\le 2\sigma_2(x),
\end{equation}
\begin{equation}\label{Convexp5}
\hspace{-5.5cm} x \sigma_2'(y)\le \sigma_2(x)+\sigma_2(y),
\end{equation}
and
\begin{equation}\label{Convexp6}
0 \le \sigma_1(x+y)-\sigma_1(x)-\sigma_1(y)\le  2\frac{x \sigma_1(y) + y \sigma_1(x)}{(x+y)},
\end{equation}
 for all $x, y \in \mathbb{R}_{>0}$.
\end{lemma}

Now, we are in a position to state the main theorems of this paper. The first two theorems deal with the existence of solutions
to the initial value problem \eqref{GEDM}--\eqref{Initialdata} with rate kernels satisfying \eqref{SumK} and \eqref{ProductK},
respectively.
\begin{theorem}\label{TheoremGEDM}
 Consider a function $\zeta^{in}$ satisfying \eqref{ConditionIntial}. In addition, $\varphi$ satisfies \eqref{PhiIntegrable}.
If the function $\mathcal{A}$ satisfies \eqref{SumK}, \eqref{Sobolevbound} and \eqref{boundderivative}, then there exists a weak solution $\zeta$  to \eqref{GEDM} and \eqref{Initialdata} such that:
 \begin{align}\label{TheoremEquation1}
  \zeta \in \mathcal{C}([0, \infty); L^1_{0, 1} ) \cap L^{\infty}(0, T;   L_{0, 1}^1 ),  \ \text{for each}\ T>0.
 \end{align}
 \end{theorem}

\begin{theorem}\label{TheoremGEDM2}
Changing only \eqref{SumK} to \eqref{ProductK} in Theorem~\ref{TheoremGEDM} the same result holds.
 \end{theorem}

We also prove that under a more restrictive condition on the rate coefficients (see \eqref{UniqueThm}) the initial value problem
\eqref{GEDM}--\eqref{Initialdata}  has unique solutions.

\begin{theorem}\label{UniquenessTheorem}
Let $\zeta$ be a weak solution to \eqref{GEDM}--\eqref{Initialdata}. If the coagulation equation satisfies the condition:
\begin{align}\label{UniqueThm}
\mathcal{A}(x, y; z) \le A (1+x)^{1/2}(1+y)^{1/2} \varphi(z),\ \   \text{for all} \ \  (x, y, z) \in \mathbb{R}_{>0}^3,
\end{align}
for some $A \ge 0$, and the function $\varphi$ is as defined in Theorem \ref{TheoremGEDM}.  Then \eqref{GEDM}--\eqref{Initialdata} admits a unique weak solution.
 \end{theorem}

The last two results, to be proved in section~5, establish that, for rate kernels with a bound of the sum type \eqref{SumK},
 the total number of clusters and the total mass are conserved quantities.

\begin{theorem}\label{TPTheorem}
  Assume $\mathcal{A}$ satisfies \eqref{SumK}. Then, every weak solution $\zeta$ to \eqref{GEDM}--\eqref{Initialdata} 
conserves the total amount of clusters initially in the system, i.e., it satisfies, for all $t>0,$
\begin{align*}
  \mathcal{M}_0(t)  =   \mathcal{M}_0^{in} :=  \int_0^{\infty}     \zeta^{in}(x)   dx.
\end{align*} 
\end{theorem}

\begin{theorem}\label{MassTheorem}
With the assumptions of Theorem~\ref{TPTheorem} every weak solution $\zeta$ to \eqref{GEDM}--\eqref{Initialdata} 
conserves mass, i.e., for all $t>0$,
\begin{align*}
 \mathcal{M}_1(t)  =   \mathcal{M}_1^{in} :=  \int_0^{\infty}   x  \zeta^{in}(x)   dx.
\end{align*} 
\end{theorem}


\section{Existence}\label{sec:exist}
In this section, our objective is to establish the proof for theorems \ref{TheoremGEDM} and \ref{TheoremGEDM2}. The strategy to accomplish these theorems is carried out by using a weak $L^1$ compactness technique. In order to apply this weak compactness technique our primary objective is to obtain a unique classical solution for each $n \in \mathbb{N}$ to the following truncated equations:

\begin{align}\label{TCGEDM}
\frac{\partial \zeta(t, x) }{\partial t}  =  \Bigl(\mathcal{B}_1^n ( \zeta(t, \cdot) )(x)  + \mathcal{D}_1^n ( \zeta(t, \cdot) )(x)\Bigr)  + 
\Bigl(\mathcal{B}_2^n ( \zeta(t, \cdot) )(x) + \mathcal{D}_2^n ( \zeta(t, \cdot) )(x)\Bigr),
\end{align}
where
\begin{align*}
 \mathcal{B}_1^n( \zeta(t, \cdot ))(x)  := & \int_0^{n-x}  \int_0^{n-z} \mathcal{A}_n(x+z, y; z)  \  \zeta(t, x+z) \zeta(t, y) dy dz, \\
 \mathcal{D}_1^n( \zeta(t, \cdot) )(x)  := & - \int_0^{n-x}  \int_z^{n} \mathcal{A}_n(y, x; z)  \  \zeta(t, y) \zeta(t, x) dy dz,   \\
 \mathcal{B}_2^n( \zeta(t, \cdot) )(x)  := & \int_0^{x} \int_z^{n} \mathcal{A}_n(y, x-z; z)  \  \zeta(t, y) \zeta(t, x-z) dy dz,  \\
 \mathcal{D}_2^n( \zeta(t, \cdot) )(x)  := & - \int_0^{x} \int_0^{n-z} \mathcal{A}_n(x, y; z)  \   \zeta(t, x) \zeta(t, y) dy dz,  
\end{align*}
where 
\begin{equation*}
\mathcal{A}_n (x, y; z):= \mathcal{A}(x, y; z) \chi_{(0, n)}(x) \chi_{(0, n)}(y+z)
\end{equation*}
and $\chi_{S}$ denotes the characteristic function of $S$, defined as:
\begin{equation*}
\chi_{S} (x):=\begin{cases}
1,\       &  \text{if}\ x \in  S,\ \\
0,\       &  \text{if}\ x \notin S,
\end{cases}   
\end{equation*}
and the initial condition is truncated as follows:
\begin{align}\label{TInitial}
\zeta_n(0, x)= \zeta_n^{in}(x) :=  \zeta^{in}(x)\chi_{(0, n)}.
\end{align}
We shall also use the notation $\varphi_n(z):=\varphi(z)\chi_{(0,n)}(z)$.

\medskip

In order to study \eqref{TCGEDM}--\eqref{TInitial} we introduce the following lemma, which will be useful to derive crucial results later on.

\subsection{Truncated Moment Equation}\label{sec:truncmom}
Let us start by  introducing the following lemma.
\begin{lemma}\label{TruncatedLemma}
 Let $\zeta_n(t, x) $  be a solution to \eqref{TCGEDM}--\eqref{TInitial} and $\omega \in L^{\infty}(0, n) $. Then we have:
\begin{align*} 
\frac{d}{dt} \int_0^{n} \omega(x) \zeta_n(t, x) dx =  \int_0^{n} \int_0^{x} \int_0^{n-z}    \tilde{\omega}(x, y, z)  \mathcal{B}_n(x, y; z, \zeta(t) ) dy dz dx,
\end{align*} 
with
\begin{align}
  \tilde{\omega}(x, y, z) :=  \omega(y+z) + \omega(x-z) -\omega(x) - \omega(y), \label{omegatil}
\end{align}
and
\begin{align*}
 \mathcal{B}_n(x, y; z, \zeta(t )) := \mathcal{A}_n(x, y; z)   \zeta_n(t, x) \zeta_n(t, y).   
\end{align*}
\end{lemma}
\begin{proof}
Let $\zeta_n $ be a solution to \eqref{TCGEDM}--\eqref{TInitial},   and $\omega \in L^{\infty}(0, n)$. First, we simplify the following term using Fubini's theorem, the transformations $x+z = x'$ and $z=z'$, and dropping the prime as follows:
\begin{align}\label{TIdentity1}
\int_0^{n}  \omega(x)  \mathcal{B}_1^n( \zeta(t, \cdot) )(x)  dx = & \int_0^{n}  \int_0^{n-x}  \int_0^{n-z} \omega(x) \mathcal{B}_n(x+z, y; z, \zeta(t) ) \  dy dz dx 
\nonumber\\
= & \int_0^{n}  \int_0^{x}  \int_0^{n-z} \omega(x-z) \mathcal{B}_n(x, y; z, \zeta(t) )   \  dy dz dx. 
\end{align}
Similarly, we can write the next term by changing notation $x \leftrightarrow y$ and applying Fubini's theorem thrice to get:
\begin{align}\label{TIdentity2}- 
\int_0^{n}  \omega(x)  \mathcal{D}_1^n( \zeta(t, \cdot) )(x)  dx = & \int_0^{n} \int_0^{n-x}  \int_z^{n} \omega(x) \mathcal{B}_n(y, x; z, \zeta(t) )   \  dy dz dx \nonumber\\
= & \int_0^{n} \int_0^{n-y}  \int_z^{n} \omega(y) \mathcal{B}_n(x, y; z, \zeta(t) )   \  dx dz dy \nonumber\\
= & \int_0^{n}  \int_0^{x} \int_0^{n-z} \omega(y) \mathcal{B}_n(x, y; z, \zeta(t) )  \ dy dz dx.
\end{align}
Next, using Fubini's theorem, the transformation $x-z=x'$ and $z=z'$ (and dropping the primes) and applying Fubini's again, we obtain:
\begin{align}\label{TIdentity3}
\int_0^{n}  \omega(x)  \mathcal{B}_2^n ( \zeta(t, \cdot) )(x)  dx = & \int_0^{n}  \int_0^{x} \int_z^{n} \omega(x) \  \mathcal{B}_n(y, x-z; z, \zeta(t) )   \ dy dz dx \nonumber\\
= & \int_0^{n} \int_z^{n} \int_z^{n}   \omega(x)    \mathcal{B}_n(y, x-z; z, \zeta(t))   \ dy dx dz \nonumber\\
= & \int_0^{n} \int_z^{n} \int_0^{n-z} \omega(y+z)  \mathcal{B}_n(x, y; z, \zeta(t) )   \ dy dx dz \nonumber\\
= & \int_0^{n} \int_0^{x} \int_0^{n-z} \omega(y+z)  \mathcal{B}_n(x, y; z, \zeta(t) )   \ dy dz dx.
\end{align}
Finally, again using Fubini's theorem we obtain:
\begin{align}\label{TIdentity4}
- \int_0^{n}  \omega(x) \ \mathcal{D}_2^n ( \zeta(t, x) )  dx = & \int_0^{n} \int_0^{x} \int_0^{n-z} \omega(x) \mathcal{B}_n(x, y; z, \zeta(t) ) \  dy dz dx.
\end{align}
Combining  \eqref{TIdentity1}---\eqref{TIdentity4}, we get:
\begin{align*}
\frac{d}{dt} \int_0^{n} \omega(x) \ \zeta_n(t, x) dx =  \int_0^{n} \int_0^{x} \int_0^{n-z}    \tilde{\omega}(x, y, z) \mathcal{B}_n(x, y; z, \zeta(t) ) \ dy dz dx.
\end{align*} 
This concludes the proof of Lemma \ref{TruncatedLemma}.
\end{proof}


\begin{proposition}\label{Prop1}
For any $n > 1$ then equation \eqref{TCGEDM}--\eqref{TInitial} has a unique non-negative solution $\zeta_n\in \mathcal{C}^1([0,\infty);L^1(0,n))$. Furthermore, this solution satisfies
\begin{align}\label{PropMassbound}
\int_0^n x \zeta_n(t, x) dx = \int_0^n x \zeta_n^{in}(x) dx, 
\end{align}
and
\begin{align}\label{PropTotalParticlebound}
\int_0^n  \zeta_n(t, x) dx = \int_0^n  \zeta_n^{in}(x) dx, 
\end{align}
for $t\ge 0$.
\end{proposition}

\begin{proof}
The strategy to prove this theorem relies on the Picard-Lindel\"{o}f theorem  \cite[Theorem 7.3]{Brezis:2011} in the Banach 
space $L^1(0, n)$. To begin, the equation \eqref{EstimateCoagRate} allows us to bound the coagulation kernel as follows:
 \begin{align}\label{bound for kernel}
\mathcal{A}_n(x, y; z) \le 4 A\ n^2 \varphi_n(z) ,\ \ \text{for }\ \ n > 1.
 \end{align}
Next, we show that each term on the right-hand side of \eqref{TCGEDM} is locally Lipschitz continuous in the space $L^1(0, n)$. For this purpose, let $\zeta$ and $\eta \in L^1(0, n) $. We then use the first term from equation \eqref{TCGEDM}, apply Fubini's theorem, the bound \eqref{bound for kernel}, and employ the transformations $x+z =x'$ and $z=z'$, to evaluate
\begin{align}\label{Lipschitz1}
\| \mathcal{B}_1^n(\zeta) - &  \mathcal{B}_1^n(\eta)  \|_{L^1(0, n)}  \nonumber\\
\le  & \int_0^{n}  \int_0^{n-x}  \int_0^{n-z}  \mathcal{A}_n(x+z, y; z) \ |  \zeta(t, x+z) \zeta(t, y) - \eta(t, x+z) \eta(t, y) | dy dz dx \nonumber\\
= & \int_0^{n}  \int_0^{n-z}  \int_0^{n-z}  \mathcal{A}_n(x+z, y; z) \  |  \zeta(t, x+z) \zeta(t, y) - \eta(t, x+z) \eta(t, y) | 
  dy dx dz \nonumber\\
=  & \int_0^{n}  \int_z^{n}  \int_0^{n-z}  \mathcal{A}_n(x, y; z) \   |  \zeta(t, x) \zeta(t, y) - \eta(t, x) \eta(t, y) |  dy dx dz \nonumber\\
\le  & 4 A\ n^2  \| \varphi \|_{L^1(0, n)} (\| \zeta\|_{L^1(0, n)}+ \|\eta\|_{L^1(0, n)}) \  \|\zeta - \eta\|_{L^1(0, n)}.
\end{align}
Similarly, using the Fubini theorem and \eqref{bound for kernel}, we estimate
\begin{align}\label{Lipschitz2}
  \| \mathcal{D}_1^n(\zeta) - &  \mathcal{D}_1^n(\eta)  \|_{L^1(0, n)}  \nonumber\\
  \le  &   \int_0^{n}  \int_z^{n} \int_0^{n-z}   \mathcal{A}_n(x, y; z)  \  |  \zeta(t, x) \zeta(t, y) - \eta(t, x) \eta(t, y) | \ dy dx dz  \nonumber\\
\le  & 4 A\ n^2  \| \varphi \|_{L^1(0, n)} (\| \zeta \|_{L^1(0, n)}+ \|\eta\|_{L^1(0, n)}) \  \|\zeta - \eta\|_{L^1(0, n)}.
\end{align}
Next, we verify the Lipschitz continuity of the following term using the Fubini theorem and applying the transformations $x-z=x'$ and $z=z'$:
\begin{align}\label{Lipschitz3}
  \| \mathcal{B}_2^n(\zeta) - &  \mathcal{B}_2^n(\eta)  \|_{L^1(0, n)}  \nonumber\\
  \le  & \int_0^{n} \int_z^{n} \int_0^{n-z}   \mathcal{A}_n(x, y; z) \  | \zeta(t, x) \zeta(t, y) - \eta(t, x) \eta(t, y) | dy dx dz \nonumber\\
\le  & 4 A\ n^2  \| \varphi \|_{L^1(0, n)} \ (\| \zeta\|_{L^1(0, n)}+ \|\eta\|_{L^1(0, n)}) \  \|\zeta - \eta\|_{L^1(0, n)}.
\end{align}
Finally, we estimate:
\begin{align}\label{Lipschitz4}
  \| \mathcal{D}_2^n(\zeta) - &  \mathcal{D}_2^n(\eta)  \|_{L^1(0, n)}  \nonumber\\
  \le  &   \int_0^{n} \int_z^{n} \int_0^{n-z}    \mathcal{A}_n(x, y; z)  \ | \zeta(t, x) \zeta(t, y) - \eta(t, x) \eta(t, y) |  dy dx dz \nonumber\\
\le  & 4 A\ n^2  \| \varphi \|_{L^1(0, n)} (\| \zeta\|_{L^1(0, n)}+ \|\eta\|_{L^1(0, n)}) \ \|\zeta - \eta\|_{L^1(0, n)}.
\end{align}
Using the estimates from \eqref{Lipschitz1}---\eqref{Lipschitz4}, it is clear that $\mathcal{B}_{1}^n(\zeta_n)$, $\mathcal{D}_{1}^n(\zeta_n)$, $\mathcal{B}_{2}^n(\zeta_n)$, and $\mathcal{D}_2^n(\zeta_n)$ are locally Lipschitz continuous functions in the space $L^1(0, n)$. 
Therefore, according to the Picard-Lindel\"{o}f theorem there exists a unique solution $\zeta_n \in \mathcal{C}^1([0, \mathcal{T} ) ; L^1(0, n) )$ to the equations \eqref{TCGEDM}--\eqref{TInitial}. Moreover, this solution is defined on a maximal interval $t \in [0, \mathcal{T} )$, where $\mathcal{T} \in (0, \infty]$.\\
The solution either persists indefinitely ($\mathcal{T} = \infty$), or it exists up to a finite time $\mathcal{T}$, with the possibility of satisfying the condition:
\begin{align}\label{norminfinity}
\lim_{t \to \mathcal{T}<\infty  } \| \zeta_n(t)\|_{L^1(0, n) } =\infty.
\end{align}

We now study the positivity of solutions. Since $\mathcal{B}_{2}^n(\zeta)$ is locally Lipschitz continuous function in $L^1(0, n)$. Thus, the positive part, $[ \mathcal{B}_{2}^n (\zeta) ]_{+},$ 
is also locally Lipschitz continuous in $L^1(0, n)$. Again by the Picard-Lindel\"{o}f theorem we can conclude that the initial value problem 
\begin{align}\label{TCSTpositive}
\frac{\partial \zeta}{\partial t}  =  \Bigl(\mathcal{B}_1^n( \zeta) + \mathcal{D}_1^n(\zeta)\Bigr)  + \Bigl([ \mathcal{B}_2^n( \zeta)]_{+}  +  \mathcal{D}_2^n(\zeta)\Bigr),
\end{align}
with the same initial data given in \eqref{TInitial} also has a unique solution.

Let $\sign_+(r)=1$, for $r \ge 0$ and $\sign_+(r)=0$, for $r < 0$. Then, from \eqref{bound for kernel}, \eqref{TCSTpositive} and $ (-\zeta)_{+} = -\zeta \sign_{+} (-\zeta) $, we infer that
\begin{align}\label{Truncated1*}
& \frac{d}{dt} \| (-\zeta_n)_{+} (t)\|_{L^1(0, n) }  \nonumber\\
 = & -\int_0^n \!\!\!\sign_{+}(-\zeta_n)(t, x) \Bigl( \mathcal{B}_1^n(\zeta ) +  \mathcal{D}_1^n(\zeta)  +  [ \mathcal{B}_2^n(\zeta) ]_{+} + \mathcal{D}_2^n(\zeta)\Bigr) dx \nonumber\\
\le   & -\int_0^n \sign_{+}(-\zeta_n)(t, x) \ \Bigl( \mathcal{B}_1^n(\zeta) +  \mathcal{D}_1^n(\zeta) + \mathcal{D}_2^n(\zeta) \Bigr) dx \nonumber\\
=: & \;I_1 + I_2 + I_3.
\end{align}
and now, after applying Fubini's theorem, a change of variables, and the assumption on the rate coefficients,  we get, for $I_1$,
\begin{align}
I_1  = & -  \int_0^n  \int_z^{n}  \int_0^{n-z} \sign_{+}(-\zeta_n)(t, x-z)  \mathcal{A}_n(x, y; z)   \zeta_n(t, x) \zeta_n(t, y) dy dx dz  \nonumber\\
\le & 4 A n^2   \int_0^n  \int_z^{n}  \int_0^{n-z}  [ -\zeta_n)(t, x) ]_{+}  \varphi(z)   \zeta_n(t, y) dy dx dz ; \nonumber
\end{align}
for $I_2$ we change notation $x\leftrightarrow y$ and apply the bound on the rate coefficient to obtain
\begin{align*}
I_2 \le 4 A n^2 \int_0^n \int_0^{n-y}  \int_z^{n}     [-\zeta_n(t, y)]_{+}  \varphi(z)  \zeta_n(t, x)  dx dz dy.
\end{align*}
Finally, for $I_3$ the bound on the rate coefficients gives
\begin{align*}
I_3 \le  4 A n^2  \int_0^n \int_0^{x} \int_0^{n-z} [-\zeta_n(t, x)]_{+}  \varphi(z) \zeta_n(t, y) dy dz dx. 
\end{align*}
Thus, plugging these bounds into \eqref{Truncated1*} we conclude
\begin{align}\label{Truncated1}
\frac{d}{dt} \| (-\zeta_n)_{+} (t)\|_{L^1(0, n) } \le 12 A n^2 \| \varphi \|_{L^1(0, n)}  \| \zeta_n(t)\|_{L^1(0, n) } \| (-\zeta_n)_{+}(t) \|_{L^1(0, n)}.
\end{align}
Solving this differential inequality we obtain
\begin{align}\label{Truncated2}
 \| (-\zeta_n)_{+} (t)\|_{L^1(0, n )} \le & \bigg(  \| (-\zeta_n^{in} )_{+} \|_{L^1(0, n)}      \bigg) \nonumber\\
 &   \ \  \ \times \exp \bigg( 12 A n^2 \| \varphi \|_{L^1(0, n)}   \int_0^t  \| \zeta_n(s)\|_{L^1(0, n)} ds \bigg).  
\end{align}
Using the non-negativity of \eqref{TInitial} and the positive part of $(-\zeta_n^{in})_{+}$ into \eqref{Truncated2}, we get
\begin{align}\label{Truncated3}
 \| (-\zeta_n)_{+} (t)\|_{L^1(0, n) } \le 0.
\end{align}
From \eqref{Truncated3}, we can deduce that $ \zeta_n(t, \cdot) \ge 0$ for all $t \in [0, \mathcal{T})$. As a result, both \eqref{TCSTpositive} and \eqref{TCGEDM} are equivalent equations. Additionally, $\zeta_n$ satisfies \eqref{PropMassbound} and \eqref{PropTotalParticlebound}, which can be easily demonstrated using Lemma \ref{TruncatedLemma}. From the non-negativity of $\zeta_n(t, \cdot)$ for all $t \in [0, \mathcal{T})$ 
we conclude that \eqref{PropTotalParticlebound} is the same as
  \begin{align}\label{Truncated6}
 \| \zeta_n( t) \|_{L^1(0, n)} = \| \zeta_n^{in} \|_{L^1(0, n)}, \ \ \forall \ t\in [0, \mathcal{T}).
\end{align}
From this identity one can see that if $ \mathcal{T} < \infty$, then
\begin{align*}
 \lim_{t \to \mathcal{T} } \| \zeta_n(t)\|_{L^1(0, n)} =  \| \zeta_n^{in} \|_{L^1(0, n)}   < \infty.
\end{align*}
Hence, \eqref{norminfinity} is not satisfied, implying that $\mathcal{T}= \infty$. This concludes the proof of Proposition \ref{Prop1}.
\end{proof}

\subsection{Weak compactness in $L^1$} \label{sec:weekcomp}

In this section, our aim is to prove that the family of solutions $\{ \zeta_n \}_{n \in\Nb}$ is relatively compact in $\mathcal{C}([0, T]; L_{0, 1}^1 )$. To achieve this we employ a weak $L^1$ compactness technique inspired by Stewart's pioneering work on continuous coagulation-fragmentation equations, \cite{Stewart:1989}. We establish that the family of solutions $\{ \zeta_n\},$ ${n > 1}$ is uniformly bounded in $L_{0, 1}^1$, which can be deduced from \eqref{PropMassbound} and \eqref{PropTotalParticlebound}:
\begin{align}\label{Uboundlemma}
\int_0^{n} (1+ x)\ \zeta_n(t, x) dx \le  \mathcal{M}_{0} (0) + \mathcal{M}_{1} (0) =: \Gamma, \ \ \text{for all}\ t \in [0, T].
\end{align}

We proceed to analyze the behavior of the truncated solution $\zeta_n$ as the particle size $x$ becomes large. This analysis is done under the sum kernel condition \eqref{SumK}. 
\begin{lemma}\label{LargeLemma}
Let $T>0$ be fixed. Consider $\zeta^{in}$ satisfying \eqref{ConditionIntial}, and assume that the coagulation rate $\mathcal{A}$ satisfies the sum kernel bound in \eqref{SumK}. Suppose also that $\varphi$ satisfies \eqref{PhiIntegrable}. Then, for every $n>1$, 
\begin{equation}\label{C(T)}
\sup_{t\in [0, T]}\int_0^{n} \sigma_1 (x) \zeta_n(t, x)\ dx \le \Xi(T),
\end{equation}
where $\Xi(T)$ is a positive constant dependent on $T$. Here, $\sigma_1$ is a convex function that complies with \eqref{Convexp2}--\eqref{Convexp6}. 
\end{lemma}

\begin{proof} Applying Lemma~\ref{TruncatedLemma} with $\omega (x) :=\sigma_1(x) \chi_{(0, n)}(x)$ and inserting we obtain
\begin{align}\label{Large1}
\int_0^n \sigma_1(x) \zeta_n(t, x) dx = & \int_0^n \sigma_1(x) \zeta_n^{in}(x)dx  \nonumber\\
& +  \int_0^t \int_0^{n} \int_z^{n} \int_0^{n-z}    \tilde{\sigma}_1(x, y, z)  \mathcal{B}_n (x, y; z, \zeta(s) )   dy dx dz ds,
\end{align}
where
   \( \tilde{\sigma}_1  (x, y, z) =  \sigma_1 (x-z) + \sigma_1 (y+z) - \sigma_1 (x) - \sigma_1 (y). \)
Since $\sigma_1$ is a monotone increasing convex function, we have
\begin{align}\label{Large2}
 \tilde{\sigma_1 } (x, y, z)  \le & \sigma_1 (x) + \sigma_1 (y+z) - \sigma_1 (x) - \sigma_1 (y) \nonumber\\
 \le & \sigma_1 (z) + \sigma_1 (y+z) - \sigma_1 (y) - \sigma_1 (z) \nonumber\\
\le & \sigma_1 (z) + 2 \frac{ y \sigma_1 (z) + z\sigma_1 (y) } {y+z}.
 \end{align}
Using \eqref{Large2} and \eqref{Convexp2} into \eqref{Large1}, we get
\begin{align}\label{Large3}
\int_0^n \sigma_1(x) \zeta_n(t, x) dx \le  & \Gamma_1 
 +  \int_0^t \int_0^{n} \int_z^{n} \int_0^{n-z}   \sigma_1 (z)   \mathcal{B}_n (x, y; z, \zeta(s) ) dy dx dz ds     \nonumber\\
 & \quad + 2 \int_0^t \int_0^{n} \int_z^{n} \int_0^{n-z}    \frac{ y \sigma_1 (z) + z \sigma_1 (y) } {y+z} \times \nonumber\\
 & \qquad\qquad \qquad \times  \mathcal{B}_n (x, y; z, \zeta(s) ) dy dx dz ds.
 \end{align}
Let us estimate the first integral on the right-hand side of \eqref{Large3}, by using \eqref{EstimateCoagRate}, \eqref{Convexp3} and \eqref{Uboundlemma} by
 \begin{align}\label{Large4}
& \int_0^t \int_0^{n} \int_z^{n} \int_0^{n-z}   \sigma_1 (z)  \mathcal{B}_n (x, y; z, \zeta(s) ) dy dx dz ds \nonumber\\
\le & A \int_0^t \int_0^{n} \int_z^{n} \int_0^{n-z}  \sigma_1 (z)   (1+x) (1+y) \varphi(z)   \zeta_n(s, x) \zeta_n(s, y) dy dx dz ds \nonumber\\
\le & A \Gamma^2 \int_0^t \int_0^{n}   \sigma_1 (z)    \varphi(z)    dz ds  \le  A \Gamma^2 \Gamma_3 t.
\end{align}
Similarly, the last integral on the right-hand side of \eqref{Large3} can be evaluated using \eqref{SumK}, \eqref{Convexp3}, and \eqref{Uboundlemma} 
\begin{align}\label{Large5}
& 2 \int_0^t \int_0^{n} \int_z^{n} \int_0^{n-z}    \frac{ y \sigma_1 (z) + z \sigma_1 (y) } {y+z} \mathcal{B}_n (x, y; z, \zeta(s) )dy dx dz ds  \nonumber\\
\le  & 2 A \int_0^t \int_0^{n} \int_z^{n} \int_0^{n-z}   \chi_{(0, 1)}(x+y)   \frac{ y \sigma_1 (z) + z \sigma_1 (y) } {y+z}    \varphi(z)   \zeta_n(s, x) \zeta_n(s, y) dy dx dz ds \nonumber\\
& + 2 A \int_0^t \int_0^{n} \int_z^{n} \int_0^{n-z}   \chi_{(1, \infty )}(x+y)   \frac{ y \sigma_1 (z) + z \sigma_1 (y) } {y+z}    (x+y) \varphi(z)   \zeta_n(s, x) \zeta_n(s, y) dy dx dz ds \nonumber\\
\le  & 2 A \Gamma^2 \int_0^t \int_0^{n}     [  \sigma_1 (z) +  \sigma_1 (1) ]  \   \varphi(z)    dz ds \nonumber\\
& + 2 A \int_0^t \int_0^{n} \int_z^{n} \int_0^{n-z}     \frac{ y \sigma_1 (z) + z \sigma_1 (y) } {y+z}    x \varphi(z)   \zeta_n(s, x) \zeta_n(s, y) dy dx dz ds \nonumber\\
& + 2 A \int_0^t \int_0^{n} \int_z^{n} \int_0^{n-z}    \frac{ y \sigma_1 (z) + z \sigma_1 (y) } {y+z}   y \varphi(z)   \zeta_n(s, x) \zeta_n(s, y) dy dx dz ds \nonumber\\
\le  & 2 A \Gamma^2     \bigl[ 3  \Gamma_3  +  \sigma_1 (1) \| \varphi \|_{L^1}  \bigr] t  +  4A \Gamma \| \varphi \|_{L^1}  \int_0^t \int_0^{n}  \sigma_1 (x) \zeta_n(s, x) dx ds.
\end{align}

\medskip

Using \eqref{Large4} and \eqref{Large5} into \eqref{Large3} and then applying the Gronwall inequality one gets
\begin{align*}
\int_0^{n}  \sigma_1 (x)  \zeta_n(t, x) dx \le  \Xi(T),
\end{align*}
where $\Xi(T) :=\bigg(  \Gamma_1  +  A \Gamma^2 [ 7   \Gamma_3 + 2 \sigma_1(1) \| \varphi \|_{L^1} ] T \bigg) e^{ 4 A \Gamma  \| \varphi \|_{L^1}  T} $,
and this completes the proof of Lemma \ref{LargeLemma}.
\end{proof}

In a similar fashion, we next proceed to analyze the behavior of the truncated solution $\zeta_n$ as the particle size $x$ becomes large for the product kernel condition \eqref{ProductK}. 
\begin{lemma}\label{LargeLemmaP}
Let $T>0$ be fixed. Consider $\zeta^{in}$ satisfying \eqref{ConditionIntial}, and assume that the coagulation rate $\mathcal{A}$ satisfies \eqref{ProductK} and \eqref{RateCoag}. Suppose also that $\varphi$ satisfies \eqref{PhiIntegrable}. Then, for every $n>1$, 
\begin{equation}\label{C(T)1}
\sup_{t\in [0, T]}\int_0^{n} \sigma_1 (x) \zeta_n(t, x)\ dx \le \Lambda(T),
\end{equation}
where $\Lambda(T)$ is a positive constant dependent on $T$ and $\sigma_1$ is a convex function that complies with \eqref{Convexp2}--\eqref{Convexp6}. 
\end{lemma}
\begin{proof} Similar to the proof of Lemma \ref{LargeLemma}, we define $\omega (x) :=\sigma_1(x) \chi_{(0, n)}(x)$ and insert it into Lemma~\ref{TruncatedLemma}. Utilizing \eqref{Large2} and \eqref{Convexp2}, we obtain
\begin{align}\label{PLarge1}
\int_0^n \sigma_1(x) \zeta_n(t, x) dx \le  & \Gamma_1 
 +  \int_0^t \int_0^{n} \int_z^{n} \int_0^{n-z}   \sigma_1 (z)   \mathcal{B}_n (x, y; z, \zeta(s) ) dy dx dz ds     \nonumber\\
 & \quad + 2 \int_0^t \int_0^{n} \int_z^{n} \int_0^{n-z}    \frac{ y \sigma_1 (z) + z \sigma_1 (y) } {y+z} \times \nonumber\\
 & \qquad\qquad \qquad \times  \mathcal{B}_n (x, y; z, \zeta(s) ) dy dx dz ds.
 \end{align}
We evaluate the first integral on the right-hand side of \eqref{PLarge1} using \eqref{EstimateCoagRate}, \eqref{Convexp3} and \eqref{Uboundlemma} by
 \begin{align}\label{PLarge2}
 \int_0^t \int_0^{n} \int_z^{n} \int_0^{n-z}   \sigma_1 (z)  \mathcal{B}_n (x, y; z, \zeta(s) ) dy dx dz ds 
  \le  A \Gamma^2 \Gamma_3 t.
\end{align}
Similarly, the last integral on the right-hand side of \eqref{PLarge1} can be evaluated using \eqref{ProductK}
\begin{align}\label{PLarge3}
& 2 \int_0^t \int_0^{n} \int_z^{n} \int_0^{n-z}    \frac{ y \sigma_1 (z) + z \sigma_1 (y) } {y+z} \mathcal{B}_n (x, y; z, \zeta(s) )dy dx dz ds  \nonumber\\
\le  & 2 A \int_0^t \int_0^{1} \int_z^{1} \int_0^{1}  \frac{ y \sigma_1 (z) + z \sigma_1 (y) } {y+z}    \varphi(z)   \zeta_n(s, x) \zeta_n(s, y) dy dx dz ds \nonumber\\
& + 2 A \int_0^t \int_0^{1} \int_z^{1} \int_1^{n}  \frac{ y \sigma_1 (z) + z \sigma_1 (y) } {y+z}     \eta(y) \varphi(z)   \zeta_n(s, x) \zeta_n(s, y) dy dx dz ds \nonumber\\
& + 2 A \int_0^t \int_0^{n} \int_1^{n} \int_0^{1}   \frac{ y \sigma_1 (z) + z \sigma_1 (y) } {y+z} \eta(x) \varphi(z)   \zeta_n(s, x) \zeta_n(s, y) dy dx dz ds \nonumber\\
& + 2 A \int_0^t \int_0^{n} \int_1^{n} \int_1^{n}     \frac{ y \sigma_1 (z) + z \sigma_1 (y) } {y+z} \eta(x) \eta(y) \varphi(z)   \zeta_n(s, x) \zeta_n(s, y) dy dx dz ds \nonumber\\
 = : & J_1 + J_2 + J_3 + J_4.
\end{align}
Using the monotonicity of $\sigma_1$, \eqref{PhiIntegrable} and \eqref{Uboundlemma}, we can write
\begin{align}\label{PLarge4}
 J_1 \leq &  4 A \sigma_1(1) \int_0^t \int_0^{1} \int_z^{1} \int_0^{1}     \varphi(z)   \zeta_n(s, x) \zeta_n(s, y) dy dx dz ds \nonumber\\
\leq & 4 A \sigma_1(1)  \Gamma^2    \| \varphi \|_{L^1} T.
\end{align}
Similarly, applying the monotonicity of $\sigma_1$, \eqref{PhiIntegrable}, \eqref{RateCoag}, and \eqref{Uboundlemma}, we get
\begin{align}\label{PLarge5}
 J_2 \leq &  2 A \int_0^t \int_0^{1} \int_z^{1} \int_1^{n}  \frac{ y \sigma_1 (1) +  \sigma_1 (y) } {y+z}     \eta(y) \varphi(z)   \zeta_n(s, x) \zeta_n(s, y) dy dx dz ds \nonumber\\
\leq & 2 A \sigma_1(1)  \Gamma^2  \eta^{\ast}  \| \varphi \|_{L^1} T + 4 A \eta^{\ast} \| \varphi \|_{L^1} \Gamma \int_0^t   \int_1^{n}     \sigma_1 (y)  \zeta_n(s, y) dy  ds.
\end{align}
We bound $ J_3$, by using the monotonicity of $\sigma_1$, \eqref{PhiIntegrable}, \eqref{RateCoag}, \eqref{Convexp3} and \eqref{Uboundlemma} by
\begin{align}\label{PLarge6}
 J_3 \leq &   2 A \int_0^t \int_0^{n} \int_1^{n} \int_0^{1}   \frac{ y \sigma_1 (z) + z \sigma_1 (y) } {y+z} \eta(x) \varphi(z)   \zeta_n(s, x) \zeta_n(s, y) dy dx dz ds \nonumber\\
\leq & 2 A  \Gamma^2 \Gamma_3 \eta^{\ast}  \| \varphi \|_{L^1} T + 2 A \sigma_1(1) \eta^{\ast} \| \varphi \|_{L^1} \Gamma^2 T.
\end{align}
Finally, $ J_4$ can be estimated, by using the monotonicity of $\sigma_1$, \eqref{PhiIntegrable}, \eqref{RateCoag}, \eqref{Convexp3} and \eqref{Uboundlemma} 
\begin{align}\label{PLarge7}
 J_4 \leq &  2 A \int_0^t \int_0^{n} \int_1^{n} \int_1^{n}     \frac{ y \sigma_1 (z) + z \sigma_1 (y) } {y+z} \eta(x) \eta(y) \varphi(z)   \zeta_n(s, x) \zeta_n(s, y) dy dx dz ds \nonumber\\
\leq & 2 A  \Gamma^2  {\eta^{\ast} }^2 \Gamma_3 T  + 4 A \Gamma^2  {\eta^{\ast} }^2 \| \varphi \|_{L^1_{0, 1}} \int_0^t  \int_0^{n} \sigma_1 (x) \zeta_n(s, x) dx ds.
\end{align}
It can be inferred from \eqref{PLarge1}--\eqref{PLarge7} that
\begin{align*}
\int_0^{n}  \sigma_1 (x) \zeta_n(t, x) dx  \le   \Theta_1(T)  +  \Theta_2(T)  \int_0^t \int_0^{n}  \sigma_1 (x) \zeta_n(s, x) dx ds,
\end{align*}
where 
\begin{align*}
 \Theta_1(T) := \Gamma_1 + 2AT  \Gamma^2  [  2 \sigma_1(1)    \| \varphi \|_{L^1}  +   2 \sigma_1(1)    \eta^{\ast}  \| \varphi \|_{L^1}   +    \Gamma_3 \eta^{\ast}  \| \varphi \|_{L^1}  +      {\eta^{\ast} }^2 \Gamma_3 ],
\end{align*}
and 
\begin{align*}
 \Theta_2(T) :=   4 A \eta^{\ast} \| \varphi \|_{L^1_{0, 1} } \Gamma (1 +  \Gamma  \eta^{\ast} ).
\end{align*}
\medskip
We then apply Gronwall's inequality to get
\begin{align*}
\int_0^{n}  \sigma_1 (x)  \zeta_n(t, x) dx \le  \Lambda(T),
\end{align*}
where $\Lambda(T) :=\bigg(   A \Gamma^2     \Gamma_3 T + \Theta_1(T)    \bigg) e^{ \Theta_2(T)  T} $.
This finishes the proof of Lemma \ref{LargeLemmaP}.
\end{proof}

The lemma in the next subsection establishes equi-integrability for the family of solutions $\{ \zeta_n \}_{n \in \Nb}$.
\subsection{Equi-integrability}\label{sec:equiint}
We now introduce the following convex function:
\begin{equation*}
 \sigma_{2, \lambda}(x) :=    \begin{cases}
\sigma_2(x), & \text{ if } x\in [0,  \lambda ],\\
\sigma_2'(\lambda) (x-\lambda) + \sigma_2(\lambda), & \text{if} \ x  \in [\lambda, \infty ).
\end{cases}
\end{equation*}
It can be readily observed that $\sigma_{2, \lambda}$ satisfies the conditions \eqref{Convexp2}--\eqref{Convexp6}.

\begin{lemma}\label{LemUformInteg}
Let $\zeta^{in}$ satisfy \eqref{ConditionIntial} and assume that the coagulation rate $\mathcal{A}$ satisfies 
the sum kernel stated in either \eqref{SumK} or \eqref{ProductK}. Suppose $\varphi(z)$ satisfies assumption \eqref{PhiIntegrable}. 
Then, for any $T>0$, there exists a positive constant $C(T)$ such that
\begin{align*}
\sup_{t\in [0,T]}\int_0^{\infty} \sigma_2( \zeta_n(t, x))dx \le C(T),
\end{align*}
where $\sigma_2$ is a convex function and $\sigma'_2$ is concave, satisfying \eqref{Convexp2}--\eqref{Convexp6}.
\end{lemma}
\begin{proof} 
Let $ \lambda >1 $. Then from \eqref{TCGEDM}, the non-negativity of $\zeta_n$, and $\mathcal{A}_n $,  we obtain
\begin{align}\label{Equint1}
&  \frac{d}{dt}\int_0^{\infty } \sigma_{2, \lambda} ( \zeta_n(t, x)) dx \nonumber\\
 = &  \int_0^{ n } \int_0^{n-x}  \int_0^{n-z}  \sigma'_{2, \lambda} ( \zeta_n(t, x))  \mathcal{B}_n (x+z, y; z, \zeta(t) )dy dz dx\nonumber\\
 & -   \int_0^{ n } \int_0^{n-x}  \int_z^{n}  \sigma'_{2, \lambda}( \zeta_n(t, x))  \mathcal{B}_n (y, x; z, \zeta(t) )dy dz dx\nonumber\\
&+\int_0^{ n }  \int_0^{x} \int_z^{n}  \sigma'_{2, \lambda}(\zeta_n(t, x)) \mathcal{B}_n (y, x-z; z, \zeta(t) )  dy dz dx\nonumber\\
& -  \int_0^{ n } \int_0^{x}  \int_0^{n-z}  \sigma'_{2, \lambda}( \zeta_n(t, x)) \mathcal{B}_n (x, y; z, \zeta(t) )  dy dz dx \nonumber\\
= : &  I_4 - I_5 + I_6 -I_7.
\end{align}
Employing Fubini's theorem and the transformation $x+z=x'$ and $z=z'$ we evaluate $I_4$ as follows
\begin{align}\label{Equint2}
I_4 =    &  \int_0^{ n} \int_0^{n-z}  \int_0^{n-z}  \sigma'_{2, \lambda}( \zeta_n(t, x)) \mathcal{B}_n (x+z, y; z, \zeta(t) ) dy dx dz \nonumber\\
=  &  \int_0^{ n } \int_z^{n}  \int_0^{n-z}  \sigma'_{2, \lambda}( \zeta_n(t, x-z))\mathcal{B}_n (x, y; z, \zeta(t) )dy dx dz. 
\end{align}
Again applying the Fubini theorem and changing the notation $x\leftrightarrow y$, we find
\begin{align}\label{Equint3}
 I_5 =    &   \int_0^{ n }   \int_0^{n-z} \int_z^{n}  \sigma'_{2, \lambda}( \zeta_n(t, x)) \mathcal{B}_n (y, x; z, \zeta(t) ) dy dx dz \nonumber\\
 =   &   \int_0^{ n }  \int_z^{n}   \int_0^{n-z}     \sigma'_{2, \lambda} ( \zeta_n(t, y)) \mathcal{B}_n (x, y; z, \zeta(t) )dy dx dz.
\end{align}
Using the Fubini theorem,  the change of variables $x\mapsto x':= x-z$ (and dropping the prime), change notation $x\leftrightarrow y$,
and finally Fubini's theorem again, we obtain
\begin{align}\label{Equint4}
I_6 =  & \int_0^{ n}  \int_z^{n} \int_z^{n}  \sigma'_{2, \lambda} (\zeta_n(t, x))  \mathcal{B}_n (y, x-z; z, \zeta(t) ) dy dx dz \nonumber\\
 =   & \int_0^{ n }   \int_0^{n-z}   \int_z^{n}  \sigma'_{2, \lambda}( \zeta_n(t, x+z)) \mathcal{B}_n (y, x; z, \zeta(t) ) dy dx dz \nonumber\\
=    & \int_0^{ n }   \int_0^{n-z}   \int_z^{n}  \sigma'_{2, \lambda}( \zeta_n(t, y+z))  \mathcal{B}_n (x, y; z, \zeta(t) )  dx dy dz\nonumber\\
=    & \int_0^{ n }   \int_z^{n}  \int_0^{n-z}    \sigma'_{2, \lambda} \zeta_n(t, y+z))  \mathcal{B}_n (x, y; z, \zeta(t) )  dy dx dz.
\end{align}
We can simplify the last term, by using only the Fubini theorem 
\begin{align}\label{Equint5}
I_7 =   &  \int_0^{ n } \int_z^{n}  \int_0^{n-z}   \sigma'_{2, \lambda}( \zeta_n(t, x))\mathcal{B}_n (x, y; z, \zeta(t) )  dy dx dz.
\end{align}
Using \eqref{Equint2}--\eqref{Equint5} into \eqref{Equint1}, we obtain
\begin{align}\label{Equint7}
&  \frac{d}{dt}\int_0^{\infty }  \sigma_{2, \lambda}( \zeta_n(t, x)) dx \nonumber\\
 \le & \int_0^{ n } \int_z^{n}  \int_0^{n-z} [   \sigma'_{2, \lambda}( \zeta_n(t, x-z)) -  \sigma'_{2, \lambda}( \zeta_n(t, x)) ] \mathcal{B}_n (x, y; z, \zeta(t) )  dy dx dz \nonumber\\
& + \int_0^{ n }   \int_z^{n}  \int_0^{n-z}   [  \sigma'_{2, \lambda}( \zeta_n(t, y+z)) -  \sigma'_{2, \lambda}( \zeta_n(t, y)) ] \mathcal{B}_n (x, y; z, \zeta(t) )  dy dx dz. 
\end{align}
From the convexity of $ \sigma_{2, \lambda}$ it is obvious that  $\sigma'_{2, \lambda}(x) \le \frac{\sigma_{2, \lambda}(x)- \sigma_{2, \lambda}(y)}{x-y}$,
which we can write as
\begin{align}
 x (  \sigma'_{2, \lambda}(y) - \sigma'_{2, \lambda}(x)  )  \le \vartheta(y) - \vartheta(x),\label{teta}
\end{align}
where
\begin{align*}
 \vartheta (x) := x  \ \sigma'_{2, \lambda}(x) -  \sigma_{2, \lambda}(x). 
\end{align*}
Using \eqref{teta} in \eqref{Equint7} and the non-negativity of $\zeta_n$ and $\mathcal{A}_n$  we get, after a  straightforward manipulation,

\begin{align}\label{Equint8}
&  \frac{d}{dt}\int_0^{\infty } \sigma_{2, \lambda}( \zeta_n(t, x)) dx \nonumber\\
 \le & \int_0^{ n } \int_z^{n}  \int_0^{n-z} [  \vartheta( \zeta_n(t, x-z)) -  \vartheta ( \zeta_n(t, x)) ] \ \mathcal{A}_n(x, y; z)    \zeta_n(t, y) dy dx dz \nonumber\\
& + \int_0^{ n }   \int_z^{n}  \int_0^{n-z}    [  \vartheta ( \zeta_n(t, y+z)) - \vartheta ( \zeta_n(t, y)) ] \  \mathcal{A}_n(x, y; z)   \zeta_n(t, x)  dy dx dz\nonumber\\
\le & \int_0^{ n } \int_0^{n}  \int_0^{n-z}   \vartheta( \zeta_n(t, x))  \ [ \mathcal{A}_n(x+z, y; z) -  \mathcal{A}_n(x, y; z) ] \  \zeta_n(t, y) dy dx dz \nonumber\\
 & + \int_0^{ n } \int_0^{z}  \int_0^{n-z}   \vartheta( \zeta_n(t, x))  \  \mathcal{A}_n(x, y; z)  \  \zeta_n(t, y) dy dx dz \nonumber\\
& + \int_0^{ n }   \int_z^{n}  \int_0^{n-z}     \vartheta ( \zeta_n(t, y))  \ [  \mathcal{A}_n(x, y-z; z) - \mathcal{A}_n(x, y ; z) ] \ \zeta_n(t, x)  dy dx dz\nonumber\\
& =: I_8 + I_9 + I_{10}.
\end{align}
Using \eqref{Sobolevbound} and \eqref{boundderivative} we can write
\begin{align*}
I_8  \le & \int_0^{ n } \int_0^{n}  \int_0^{n-z}   \vartheta( \zeta_n(t, x)) \bigg| \frac{ \partial}{\partial x} { \{ \mathcal{A}_n(x, y; z) \} } \bigg|  \  z   \zeta_n(t, y) dy dx dz \nonumber\\
\le & A  \int_0^{ n } \int_0^{n}  \int_0^{n}   \vartheta( \zeta_n(t, x))  (1+y)   z \varphi_n(z)  \zeta_n(t, y) dy dx dz \nonumber\\
\le & A \Gamma \| \varphi \|_{L^1} \int_0^{ \infty }    \vartheta( \zeta_n(t, x))   dx\nonumber \\
\le & A \Gamma \| \varphi \|_{L^1} \int_0^{ \infty } \sigma_{2, \lambda} ( \zeta_n(t, x))   dx.
\end{align*}
and the same holds true for $I_{10}$. The physical assumption \eqref{physical} immediately implies that $I_9=0.$

Thus, for any $T>0$ fixed, we can apply the Gronwall inequality and obtain:
\begin{align*}
 \int_0^{ \infty } \sigma_{2, \lambda}( \zeta_n(t, x)) dx \le \bigg( \int_0^{ \infty } \sigma_{2, \lambda} ( \zeta_n^{in}(x)) dx \bigg)   e^{ 2 A \Gamma \| \varphi \| T}.
\end{align*}
If we take the limit as $\lambda \to \infty$ the proof is completed. 

\end{proof}


\subsection{Equi-continuity with respect to time }\label{sec:timeequicont}
In order to apply the Arzel\`a-Ascoli theorem to prove existence we need to establish an equicontinuity result for the sequence $\{ \zeta_n(t, x) \}_{n >1}.$

\begin{lemma}\label{Equicontinuity}
Consider the assumptions of Lemma~\ref{LemUformInteg}.
For any $T>0$ and $ \lambda  >1$, there exists a positive constant $C_0$ such that
\begin{align*}
\int_0^{\lambda}   | \zeta_n(t, x)- \zeta_n(s, x) | dx  \le C_0 (t-s),
\end{align*}
for every $n > 1$, $0 \le s \le t \le T$.
\end{lemma}
\begin{proof}
Assume $T>0$ and $\lambda >1$. For $0 \le s \le t \le T$, we evaluate 
\begin{align}\label{Equicontinuity1}
 \int_0^{\lambda}    |\zeta_n(t, x)- \zeta_n(s, x)| dx 
 \le  &  \int_s^t  \bigg[   \int_0^{\lambda}  \mathcal{B}_1^n(\zeta (\tau, \cdot ) )(x) dx  + \int_0^{\lambda}  |\mathcal{D}_1^n(\zeta (\tau, \cdot) )(x)| dx \nonumber\\
  &+ \int_0^{\lambda}   \mathcal{B}_2^n(\zeta (\tau, \cdot) )(x) dx  + \int_0^{\lambda}   |\mathcal{D}_2^n(\zeta (\tau, \cdot) )(x)| dx   \bigg] d\tau.
\end{align}
By applying the Fubini theorem and  using the transformation $x+z=x'$ and $z=z'$ (and then dropping the primes), the first term on the right-hand side to \eqref{Equicontinuity1} can be evaluated as
\begin{align}\label{Equicontinuity2}
\int_0^{\lambda}  \mathcal{B}_1^n( \zeta(\tau, \cdot) )(x) dx 
= &  \int_0^{n-\lambda}  \int_0^{\lambda}  \int_0^{n-z} \mathcal{B}_n (x+z, y; z, \zeta(\tau) )  dy dx dz \nonumber\\
&\quad + \int_{n-\lambda}^n  \int_{0}^{n-z}  \int_0^{n-z} \mathcal{B}_n (x+z, y; z, \zeta(\tau ) )  dy dx dz \nonumber\\
= &  \int_0^{n-\lambda}  \int_z^{\lambda+z}  \int_0^{n-z} \mathcal{B}_n (x, y; z, \zeta(\tau ) )  dy dx dz \nonumber\\
&\quad + \int_{n-\lambda}^{n}  \int_{z}^{n}  \int_0^{n-z} \mathcal{B}_n (x, y; z, \zeta(\tau ) )  dy dz dx \nonumber\\
\le &  \int_0^{n-\lambda}  \int_0^{n}  \int_0^{n}  \mathcal{B}_n(x, y; z, \zeta(\tau ) )    dy dx dz \nonumber \\
&\quad + \int_{n-\lambda}^{n}  \int_{0}^{n}  \int_0^{n} \mathcal{B}_n(x, y; z, \zeta(\tau) )   dy dz dx \nonumber\\
= &  \int_0^{n}  \int_0^{n}  \int_0^{n}  \mathcal{B}_n(x, y; z, \zeta(\tau ) )    dy dx dz.
\end{align}
Then, we can estimate the above term, using \eqref{EstimateCoagRate} and \eqref{Uboundlemma}, as
\begin{align}\label{Equicontinuity3}
\int_0^{\lambda}  \mathcal{B}_1^n( \zeta(\tau, \cdot) )(x) dx \le &  A \int_0^{n}  \int_0^{n}  \int_0^{n} (1+x) (1+y) \varphi_n(z)   \zeta_n(\tau, x) \zeta_n(\tau, y) dy dx dz \nonumber\\
\le &  A \| \varphi \|_{L^1}   \Gamma^2. 
\end{align}
Next, we evaluate the second term on the right-hand side to \eqref{Equicontinuity1} by applying \eqref{EstimateCoagRate} and \eqref{Uboundlemma} as
\begin{align}\label{Equicontinuity4}
\int_0^{\lambda} |\mathcal{D}_1^n( \zeta(\tau, \cdot) )(x)| dx = & \int_0^{\lambda} \int_0^{n-x}  \int_z^{n} \mathcal{B}_n (y, x; z, \zeta(\tau ) )  dy dz dx\nonumber\\    
\le & A \int_0^{\lambda} \int_0^{n}  \int_0^{n} (1+x) (1+y) \varphi_n(z)   \zeta_n(\tau, y) \zeta_n(\tau, x) dy dz dx \nonumber\\  
\le &  A \| \varphi  \|_{L^1} \Gamma^2.
\end{align}
We can now estimate the third term  on the right-hand to \eqref{Equicontinuity1} using the Fubini theorem,  the transformation $x-z =x'$ and $z=z'$, \eqref{EstimateCoagRate} and \eqref{Uboundlemma} as
\begin{align}\label{Equicontinuity5}
\int_0^{\lambda} \mathcal{B}_2^n( \zeta(\tau, \cdot) )(x) dx = & \int_0^{\lambda}  \int_0^{x} \int_z^{n} \mathcal{B}_n (y, x-z; z, \zeta(\tau) )dy dz dx\nonumber\\    
\le  & \int_0^{\lambda}  \int_z^{\lambda} \int_0^{n} \mathcal{B}_n (y, x-z; z, \zeta(\tau) )dy dx dz\nonumber\\  
\le  & A \int_0^{\lambda}  \int_0^{\lambda} \int_0^{n}   (1+x) (1+y)  \varphi_n(z) \zeta_n(\tau, x) \zeta_n(\tau, y) dy dx dz\nonumber\\  
\le  &  A \| \varphi  \|_{L^1}  \Gamma^2. 
\end{align}
Finally, the last term on the right-hand side to \eqref{Equicontinuity1} can be evaluated by applying \eqref{EstimateCoagRate} and \eqref{Uboundlemma} as
\begin{align}\label{Equicontinuity6}
\int_0^{\lambda}|\mathcal{D}_2^n( \zeta(\tau, \cdot) )(x)| dx =  & \int_0^{\lambda} \int_0^{x} \int_0^{n-z} \mathcal{B}_n(x, y; z, \zeta(\tau) )  dy dz dx \nonumber\\
\le & A \| \varphi \|_{L^1}   \int_0^{\lambda} \int_0^{n} (1+x) (1+y)   \varphi_n(z) \zeta_n(\tau, x) \zeta_n(\tau, y) dy  dx \nonumber\\
\le &  A \| \varphi \|_{L^1}     \Gamma^2. 
\end{align}
Using \eqref{Equicontinuity3}, \eqref{Equicontinuity4}, \eqref{Equicontinuity5} and \eqref{Equicontinuity6} into \eqref{Equicontinuity1}, we get
\begin{align*}
\int_0^{\lambda}    | \zeta_n(t, x) - \zeta_n(s, x)| dx \le  C_0 (t-s),
\end{align*}
where
\begin{align*}
C_0 : = 4 A \| \varphi \|_{L^1}     \Gamma^2.
\end{align*}
This completes the proof of the Lemma \ref{Equicontinuity}.
\end{proof}

Based on Lemma \ref{Equicontinuity}, we ensure that the family $\{ \zeta_n \}_{n\in\Nb}$ exhibits strong continuity. Given that strong continuity implies weak continuity, it follows that the family $\{ \zeta_n \}_{n\in\Nb}$ is also weakly continuous in $\mathcal{C}([0, T], w\text{-}L^1(0, \lambda))$. In other words, 
for every $\omega \in L^{\infty }$, we have:\\
\begin{align}\label{weakcontinuity}
\int_0^{\lambda}  \omega(x) | \zeta_n(t, x) - \zeta_n(s, x) | dx  \le \| \omega \|_{L^{\infty }(\mathbb{R}_{>0})} C_0 (t-s).
\end{align}

We are now in a position to complete the proof of Theorem \ref{TheoremGEDM} in the next subsection.

\subsection{Convergence of integrals for the sum kernel.}\label{sec:convsum}
We are now ready to prove Theorem~\ref{TheoremGEDM}.
\begin{Proof} \emph{of Theorem}~\ref{TheoremGEDM}.
By employing a refined version of the de la Vall\'{e}e-Poussin theorem for $p=1$, \cite[Theorem~2.29]{Fonseca:2007}, 
along with Lemma \ref{LemUformInteg}, and subsequently applying the Dunford-Pettis theorem and a variant of the 
Arzel\`a-Ascoli theorem (as mentioned in \cite{Vrabie:1995}), we can deduce that the sequence $(\zeta_n)$ is relatively compact 
in $\mathcal{C}([0, T]; w\text{-}L^1(0, \lambda))$ for each $T>0$. As a result, there exists a subsequence of $(\zeta_n)$ 
(which is not relabeled) and a non-negative function $\zeta \in \mathcal{C}([0, T]; w\text{-}L^1(0, \lambda ))$ such that:
\begin{align}\label{weakconvergence}
 \zeta_n \to \zeta \ \ \ \text{in} \ \ \mathcal{C}([0, T]; w\text{-}L^1(0, \lambda)).
\end{align}
Continuing along the same lines as detailed in \cite{Barik:2020}, we can improve the convergence \eqref{weakconvergence} to
\begin{align}\label{weakconvergence1}
 \zeta_n \to \zeta \ \ \ \text{in} \ \ \mathcal{C}([0, T]; w\text{-}L_{0,1}^1 )\ \ \text{for each $T>0$},
\end{align}
by applying Lemma~\ref{LargeLemma} and \eqref{weakconvergence}.  \\
 Next, our claim is that
\begin{align}\label{Convergencestrongly}
 \zeta \in   \mathcal{C} ( [0, T]; L_{0, 1}^1 ).
\end{align}
For $t \ge s \ge 0$, applying \eqref{weakconvergence1}, we have that $\{ \zeta_n(t) - \zeta_n(s) \}$ converges weakly to 
$\{ \zeta(t) - \zeta(s) \}$ in $L^1$. This allows us to take the limit as $n \to \infty$ in \eqref{weakcontinuity}, 
resulting in $\zeta$ also satisfying \eqref{weakcontinuity}. This leads us to
 \begin{align*}
& \|\zeta(t)- \zeta(s) \|_{L^1} \nonumber\\
  = & \sup_{ \omega \in L^{\infty}  } \bigg\{  \frac{1}{\|\omega \|_{L^{\infty}   }} \bigg| \int_0^{\infty} \{ \zeta( t, x)- \zeta(s, x) \} \omega(x) dx \bigg| \bigg\} \nonumber\\
 \le  & C_0(T, \lambda) (t-s).
  \end{align*}
  This proves claim \eqref{Convergencestrongly}.  
  To conclude the proof of Theorem \ref{TheoremGEDM} it remains to prove that $\zeta$ is indeed a solution to 
\eqref{GEDM}--\eqref{Initialdata} in a weak sense. Thus, we must verify the convergence of all truncated integrals in 
\eqref{TCGEDM} to the original integrals in \eqref{GEDM}, respectively. We follow the same standard procedure to prove the integral convergence as in the Smoluchowski coagulation equation, as documented in references \cite{Barik:2020, Laurencot:2002, Laurencot:2015, Stewart:1989} and therein. 

\medskip

Suppose $\omega \in L^{\infty}$ with compact support included in $\mathbb{R}_{>0}$. 
By \eqref{TCGEDM}, we have
\begin{align}\label{Converterm1}
\int_0^{\infty} \omega(x) [\zeta_n(t, x)- \zeta_n^{in}(x)] dx  = \sum_{k=11}^{14}I_k,
\end{align}
where
$$I_k = I_k^n(t) := \iiiint_{\Omega_k}\tilde{\omega}(x, y, z) \ \mathcal{B}_n (x, y; z, \zeta(s) ) \, dz  dy dx ds,$$
with 
$\Omega_k = \bigl\{(s,x,y,z): (s,x,y)\in [0, t]\times R_k\;\text{and $z\in [0, x]$}\bigr\}$ and $R_k$ as in Figure~2,
where  $\lambda>1$ is arbitrarily fixed.

\medskip

%
%
%
\begin{figure}[h]
	\includegraphics[scale=1.3]{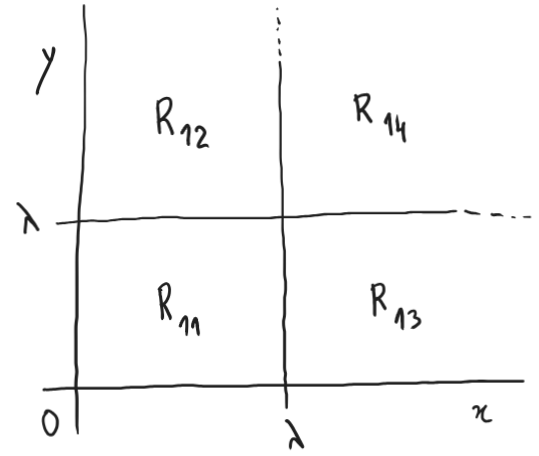}\label{figure2}
	\caption{Regions $R_k$ used in the integrals in the right-hand side of $\eqref{Converterm1}$.}
\end{figure}
%
%
%

To estimate the limit as $n\to\infty$ of the left-hand side of \eqref{Converterm1} we estimate the limit of each of the integrals $I_k$, with fixed $\lambda >1$ and then pass to the limit $\lambda\to \infty$.

\medskip

For the integral $I_{11}$  we observe that, as $n\to\infty$, $\mathcal{A}_n (x, y; z)\to  \mathcal{A}(x, y; z)$ a.e. and $\zeta_n \to \zeta $ in $\mathcal{C}([0,T]; w\text{-}L^1(\mathbb{R}_{>0}) )$, and considering condition \eqref{SumK} is 
evident that $\left\{\int_0^x \tilde{\omega}(x, y, z)  \mathcal{A}_n (x, y; z) dz\right\}$ is a bounded sequence in 
$L^{\infty} ( (0, \lambda)^2 \times (0, t) )$. Then, using \cite[Proposition 18]{Laurencot:2015} we can assert that 
\begin{multline*}
\lim_{n \to \infty}   \int_0^{\lambda }  \int_0^{x} \tilde{\omega}(x, y, z) \ \mathcal{A}_n (x, y; z)  \zeta_n(s, y) \ dz  dy  \\
 = \int_0^{\lambda }  \int_0^{x} \tilde{\omega}(x, y, z) \ \mathcal{A} (x, y; z)  \zeta(s, y) \ dz  dy,
\end{multline*}
and the same argument and the bound  $\int_0^{\lambda}  \zeta_n(t, x) dx \le \Gamma$ imply that
\begin{align}\label{J1}
\lim_{n \to \infty}I_{11} = & \lim_{n \to \infty} \int_0^t\int_{0}^{\lambda}  \zeta_n(s, x)   \bigg[   \int_{0}^{\lambda}  \int_0^x \tilde{\omega}(x, y, z) \ \mathcal{A}_n (x, y; z)  \zeta_n( s, y) dz dy  \bigg]  dx ds\nonumber \\
 = & \int_0^t \int_{0}^{\lambda} \int_{0}^{\lambda}   \int_0^x \tilde{\omega}(x, y, z) \ \mathcal{A} (x, y; z)   \zeta(s, x) \zeta(s, y) dz dy  dx ds.
\end{align}
Next, by applying \eqref{Uboundlemma}, Lemma \ref{LargeLemma}, \eqref{Convexp2} and \eqref{SumK}, $I_{12}$ can estimate as follows:
\begin{align*}
|I_{12}| = & \int_0^t \int_0^{\lambda} \int_{\lambda }^{\infty}  \int_0^{x} \tilde{\omega}(x, y, z) \ \mathcal{A}_n (x, y, z)   \ \zeta_n(s, x) \zeta_n(s, y)  dz  dy dx ds \nonumber\\
\le &  4A \|\omega\|_{L^{\infty}}  \int_0^t \int_0^{\lambda} \int_{\lambda }^{\infty}  \int_0^{x}  (x+y) \varphi_n(z)    \zeta_n(s, x) \zeta_n(s, y)  dz  dy dx ds \nonumber\\
 \le & 4A \|\omega\|_{L^{\infty}}  \| \varphi \|_{L^1} \Gamma  \int_0^t  \int_{\lambda}^{\infty}  y   \zeta_n(s, y)    dy  ds \nonumber\\
\le & 4A \|\omega\|_{L^{\infty}}  \| \varphi \|_{L^1} \Gamma   \frac{\lambda}{\sigma_1(\lambda)} \int_0^t  \int_{\lambda}^{\infty} \sigma_1(y) \zeta_n(s, y) dy  ds\nonumber\\
\le & 4A \|\omega\|_{L^{\infty}}  \| \varphi \|_{L^1} \Gamma   \frac{\lambda}{\sigma_1(\lambda)} \Xi(T)  t.
\end{align*}
Then, as $\lambda \to \infty$ we obtain
\begin{align}\label{J21}
|I_{12}| \to 0,
\end{align}
uniformly in $n$ and for $t\in [0,T]$. 

\medskip

Similarly, for $I_{13}$ we can establish that
\begin{multline*}
\bigg| \int_0^t \int_0^{\lambda} \int_{\lambda }^{\infty}  \int_0^{x} \tilde{\omega}(x, y, z) \ \mathcal{A} (x, y; z)   \ \zeta(s, x) \zeta(s, y)  dz  dy dx ds \bigg| \nonumber\\
\le  4A \|\omega\|_{L^{\infty} }  \| \varphi \|_{L^1} \| \zeta \|_{L^1_{0, 1} }   \int_0^t  \int_{\lambda}^{\infty} x \zeta(s, x) dx ds,
\end{multline*}
and, as $\lambda\to\infty$, we can apply the Lebesgue-dominated convergence theorem, leading to the conclusion that
\begin{align}\label{J22}
\int_0^t \int_0^{\lambda} \int_{\lambda }^{\infty}  \int_0^{x} \tilde{\omega}(x, y, z) \ \mathcal{A} (x, y; z)   \ \zeta(s, x) \zeta(s, y)  dz  dy dx ds  \to 0,.
\end{align}
Using \eqref{Uboundlemma}, Lemma~\ref{LargeLemma}, \eqref{SumK} and \eqref{Convexp2}, we can evaluate $I_{13}$ as follows:
\begin{align*}
|I_{13}| = & \bigg|\int_0^t \int_{\lambda}^{\infty} \int_0^{\lambda }  \int_0^{x} \tilde{\omega}(x, y, z) \ \mathcal{A}_n(x, y;  z) \zeta_n(s, x) \zeta_n(s, y) \ dz  dy dx ds \bigg| \nonumber\\
\le &  4 A \|\omega\|_{L^{\infty} }  \int_0^t \int_{\lambda}^{\infty} \int_0^{\lambda }  \int_0^{x}  (x+y) \varphi_n(z) \zeta_n(s, x) \zeta_n(s, y)   \ dz  dy dx ds \nonumber\\
\le &  4  A\Gamma  \|\omega\|_{L^{\infty} } \| \varphi \|_{L^1}    \int_0^t \int_{\lambda}^{\infty}   x \zeta_n(s, x)   \  dx ds,
\end{align*}
and as $\lambda \to \infty$ we have
\begin{align}\label{J31}
| I_{13}| \to 0.
\end{align}
Similarly to \eqref{J22}, as $ \lambda \to \infty$, the Lebesgue-dominated convergence theorem can be applied, yielding:
\begin{align}\label{J32}
\bigg|\int_0^t \int_{\lambda}^{\infty} \int_0^{\lambda }  \int_0^{x} \tilde{\omega}(x, y, z) \ \mathcal{A}(x, y;  z) \zeta(s, x) \zeta(s, y) \ dz  dy dx ds \bigg|\to 0.
\end{align}

\medskip

Finally, let us consider $I_{14}$.
\begin{align}\label{J4}
|I_{14} | \le  & 4 A \| \omega \|_{L^{\infty}}  \int_0^t \int_{\lambda}^{\infty} \int_{\lambda }^{\infty}  \int_0^{x}  (x+y)\varphi_n(z) \zeta_n(s, x) \zeta_n(y, s)    dz  dy dx ds\nonumber\\
\le  & 4 A \| \omega \|_{L^{\infty}} \| \varphi \|_{L^1} \int_0^t \int_{\lambda}^{\infty} \int_{\lambda }^{\infty}   (x+y) \zeta_n(s, x) \zeta_n(y, s)     dy dx ds\nonumber\\
\le  & 4 A \Gamma \| \omega \|_{L^{\infty}} \| \varphi \|_{L^1}  \int_0^t \int_{\lambda}^{\infty}    x \zeta_n(s, x)  dx ds.
\end{align}
As $ \lambda \to \infty$, again the Lebesgue-dominated convergence theorem can be used to establish that, as $\lambda \to \infty$,
\begin{align}\label{J41}
| I_{14} | \to 0.
\end{align}
Similarly,  we can prove that
\begin{align*}
& \int_0^t \int_{\lambda}^{\infty} \int_{\lambda }^{\infty}  \int_0^{x} \tilde{\omega}(x, y, z) \mathcal{B}_n (x, y; z, \zeta(s))   \ dz  dy dx ds \nonumber\\
\le  & 4 A \| \omega \|_{L^{\infty}}  \int_0^t \int_{\lambda}^{\infty} \int_{\lambda }^{\infty}  \int_0^{x}  (x+y)\varphi_n(z) \zeta_n(s, x) \zeta_n(y, s)    dz  dy dx ds < \infty.
\end{align*}
Thus, we have 
\begin{align}\label{J42}
\lim_{\lambda \to \infty} \bigg|  \int_0^t \int_{\lambda}^{\infty} \int_{\lambda }^{\infty}  \int_0^{x} \tilde{\omega}(x, y, z) \mathcal{B}_n(x, y; z, \zeta(s) )   \ dz  dy dx ds \bigg| =0.
\end{align}
Finally, combining \eqref{J1}, \eqref{J21}, \eqref{J31}, \eqref{J32}, \eqref{J41}, \eqref{J42} and sequentially applying the limits of $n \to \infty$ followed by $\lambda \to \infty$, we arrive at
\begin{align}\label{last}                                                        
& \lim_{n\to \infty} 
   \int_0^t \int_0^{\infty} \int_0^{\infty }  \int_0^{x} \tilde{\omega}(x, y, z) \mathcal{B}_n (x, y; z, \zeta(s) )   \ dz  dy dx ds \nonumber\\
&\qquad = \int_0^t \int_0^{\infty} \int_0^{\infty }  \int_0^{x} \tilde{\omega}(x, y, z) \mathcal{B}(x, y; z, \zeta(s) )   \ dz  dy dx ds.
\end{align}
Using once again the weak convergence, where $\zeta_n \longrightarrow \zeta$ in the space $\mathcal{C}([0,T]; w\text{-}L^1_{0, 1}( \mathbb{R}_{>0} )$, we can insert this into equation \eqref{Converterm1}. By then invoking \eqref{last}, we can deduce that
\begin{align*}
\int_0^{\infty} \omega(x) & [ \zeta(t, x)- \zeta^{in}(x) ] dx = \lim_{n \to \infty} \int_0^n \omega(x) [ \zeta_n(t, x)- \zeta_n^{in}(x)] dx \nonumber\\
 = &\lim_{n\to \infty} 
   \int_0^t \int_0^{\infty} \int_0^{\infty }  \int_0^{x} \tilde{\omega}(x, y, z) \mathcal{B}_n (x, y; z, \zeta(s) )   \ dz  dy dx ds \nonumber\\
=& \int_0^t \int_0^{\infty} \int_0^{\infty }  \int_0^{x} \tilde{\omega}(x, y, z) \mathcal{B}(x, y; z, \zeta(s) )   \ dz  dy dx ds,
\end{align*}
for every $\omega \in L^{\infty}$. This confirms that $ \zeta$ is a weak solution to \eqref{GEDM}--\eqref{Initialdata} in the sense of Definition \ref{Definition of Soln}.
\end{Proof}

\subsection{Convergence of integrals for the product kernel}\label{sec:convprod}

Our goal here is to obtain estimates (as described for the sum kernel) for the product kernel.
\begin{Proof}
\emph{of Theorem}~\ref{TheoremGEDM2}.
As in the proof of the previous theorem, what was done up to section~3.4 establishes, also for the rate kernels
satisfying \eqref{ProductK}, that there exists a sequence $(\zeta_n)$
and a nonnegative function $\zeta \in \mathcal{C}([0, T]; w\text{-}L^1((0, \lambda)))$ such that
\begin{align}\label{weakconvergencep}
 \zeta_n \to \zeta \ \ \ \text{in} \ \ \mathcal{C}([0, T]; w\text{-}L^1(0, \lambda)),
\end{align}
and 
by applying Lemma \ref{LargeLemmaP} and \eqref{Convexp1}, the above convergence can be improved to
\begin{align}\label{weakconvergencep1}
 \zeta_n \to \zeta \ \ \ \text{in} \ \ \mathcal{C}([0, T]; w\text{-}L^1_{0, 1} ). 
\end{align}
Analogously, we can establish that
\begin{align}\label{Convergencestronglyp}
 \zeta \in   \mathcal{C} ( [0, T]; L^1_{0, 1}   ).
\end{align}
We now estimate the integrals in \eqref{Converterm1} for the present rate coefficients.
The limit of $I_{11}$ can be shown to be similar to the sum kernel case, i.e. 
\begin{align}\label{PJ1}
& \lim_{n \to \infty}  \int_0^t \int_0^{\lambda} \int_{0}^{\lambda}  \int_0^{x} \tilde{\omega}(x, y, z) \ \mathcal{B}_n (x, y; z, \zeta(s) ) \ dz  dy dx ds \nonumber\\
= & \int_0^t \int_0^{\lambda}  \int_{0}^{\lambda}    \int_0^{x} \tilde{\omega}(x, y, z) \ \mathcal{B} (x, y; z, \zeta(s) ) \ dz  dy dx ds.
\end{align}
For $I_{12}$ we have also a similar estimate:
\begin{align*}
|I_{12}| = & \int_0^t \int_0^{\lambda} \int_{\lambda }^{\infty}  \int_0^{x} \tilde{\omega}(x, y, z) \ \mathcal{B}_n (x, y; z, \zeta(s) )  dz  dy dx ds \nonumber\\
\le &  4 A  \|\omega\|_{L^{\infty}}  \int_0^t \int_0^{\lambda} \int_{\lambda }^{\infty}  \int_0^{x}  \eta (x) \eta(y) \varphi(z)   \ \zeta_n(s, x) \zeta_n(s, y)  dz  dy dx ds \nonumber\\
 \le & 4 A \|\omega\|_{L^{\infty}}    \| \varphi \|_{L^1}   \int_0^t \int_0^{\lambda}  \int_{\lambda}^{\infty}  \eta (x)(1+y) \frac{ \eta(y) }{(1+y)}  \ \zeta_n(s, x) \zeta_n(s, y)    dy dx ds \nonumber\\
  \le & 4  A\Gamma \eta^{\ast} \|\omega\|_{L^{\infty}}  \| \varphi \|_{L^1}   \int_0^t \int_{\lambda}^{\infty}  (1+y)  \zeta_n(s, y)    dy  ds, 
\end{align*}
and so, as $\lambda \to \infty$ we have
\begin{align}\label{PJ21}
|I_{12}| \to 0.
\end{align}
Similarly, we can show that
\begin{align*}
& \bigg| \int_0^t \int_0^{\lambda} \int_{\lambda }^{\infty}  \int_0^{x} \tilde{\omega}(x, y, z) \ \mathcal{B} (x, y; z, \zeta(s) ) dz  dy dx ds \bigg| \nonumber\\
\le & 4A \|\omega\|_{L^{\infty} }  \| \varphi \|_{L^1}  \| \zeta \|  \int_0^t  \int_{\lambda}^{\infty}  \zeta(s, y) dy ds.
\end{align*}
As $\lambda \to \infty$, then by the Lebesgue-dominated convergence theorem, we find
\begin{align}\label{PJ22}
\lim_{\lambda \to \infty} \bigg| \int_0^t \int_0^{\lambda} \int_{\lambda }^{\infty}  \int_0^{x} \tilde{\omega}(x, y, z) \ \mathcal{A} (x, y; z)   \ \zeta(s, x) \zeta(s, y)  dz  dy dx ds \bigg|=0.
\end{align}

By using \eqref{Uboundlemma}, Lemma \eqref{LargeLemma}, \eqref{SumK} and \eqref{Convexp2}, $I_{13}$ is evaluated as
\begin{align*}
|I_{13}| = & \bigg|  \int_0^t \int_{\lambda}^{\infty} \int_{0}^{\lambda} \int_0^x \omega(x, y, z) \mathcal{B}_n (x, y; z, \zeta(s) ) dz dy dx ds \bigg| \nonumber\\
\le & 4 A \|\omega\|_{L^{\infty} } \| \varphi \|_{L^1} \int_0^t \int_{\lambda}^{\infty} \int_{\lambda}^{\infty} \eta(x) \eta(y)  \zeta_n(s, x) \zeta_n(s, y) dy dx ds \nonumber\\
\le &  4 A\eta^{\ast}\frac{\eta(\lambda)} {\lambda} \|\omega\|_{L^{\infty} }      \| \varphi \|_{L^1}   \int_0^t  \int_{\lambda}^{\infty}  (1+x)  \zeta_n(s, x)  dx ds.
\end{align*}
As $\lambda \to \infty$, we have
\begin{align}\label{PJ31}
| I_{13}| \to 0.
\end{align}
Similarly to \eqref{PJ22}, as $\lambda \to \infty$, then by the Lebesgue-dominated convergence theorem, we obtain that
\begin{align}\label{PJ32}
\bigg| \int_0^t \int_{\lambda}^{\infty} \int_{\lambda}^{\infty} \int_0^x \tilde{\omega}(x, y, z) \mathcal{B} (x, y; z, \zeta(s) ) dz dy dx ds \bigg|\to 0.
\end{align}

Let us now estimate $I_{14}$ as 
\begin{align*}
|I_{14}|  \le &  4 A  \|\omega\|_{L^{\infty}}  \int_0^t \int_{\lambda}^{\infty} \int_{\lambda }^{\infty}  \int_0^{x}  \eta (x) \eta(y) \varphi(z)   \ \zeta_n(s, x) \zeta_n(s, y)  dz  dy dx ds \nonumber\\
\le &  4 A\eta^{\ast} \|\omega\|_{L^{\infty} }    \| \varphi \|_{L^1}   \int_0^t \int_{\lambda}^{\infty}  \int_{\lambda}^{\infty}  (1+x) (1+y)  \ \zeta_n(s, x) \zeta_n(s, y)    dy dx ds \nonumber\\
  \le & 4 A\Gamma   \|\omega\|_{L^{\infty}} { \eta^{\ast} }^2  \| \varphi \|_{L^1}   \int_0^t \int_{\lambda}^{\infty}  (1+x)  \zeta_n(s, x)    dx  ds. 
\end{align*}
As $\lambda \to \infty$, then by the Lebesgue-dominated convergence theorem, we get
\begin{align}\label{PJ41}
|I_{14}| \to 0.
\end{align}
Similarly, we can show that
\begin{align*}
& \bigg| \int_0^t \int_{\lambda}^{\infty} \int_{\lambda }^{\infty}  \int_0^{x} \tilde{\omega}(x, y, z) \ \mathcal{B}(x, y; z, \zeta(s) )  dz  dy dx ds \bigg| \nonumber\\
&\qquad\qquad \le 4 \|\omega\|_{L^{\infty} (\mathbb{R}_{>0})} A \| \varphi \|_{L^1}  \| \zeta \|_{L^1_{0, 1} }  \int_0^t  \int_{\lambda}^{\infty}  \zeta(s, x) dx ds.
\end{align*}
As $\lambda \to \infty$, then by the Lebesgue dominated convergence theorem, we find
\begin{align}\label{PJ42}
\lim_{\lambda \to \infty} \bigg| \int_0^t \int_{\lambda}^{\infty} \int_{\lambda }^{\infty}  \int_0^{x} \tilde{\omega}(x, y, z)  \mathcal{B}_n(x, y; z, \zeta(s) )  dz  dy dx ds \bigg|=0.
\end{align}

Finally, combining \eqref{PJ1}, \eqref{PJ21},  \eqref{PJ22}, \eqref{PJ31}, \eqref{PJ32}, \eqref{PJ41}, \eqref{PJ42} and taking first the limit as $n \to \infty$ and then as $ \lambda \to \infty$, we end up with
\begin{align*}                                                       
& \lim_{n\to \infty}
  \int_0^t \int_0^{\infty}\int_0^{\infty} \int_0^x \tilde{\omega} (x, y, z)   \mathcal{B}_n(x, y; z, \zeta(s) ) dz dy dx ds \nonumber\\
&\qquad\qquad =  \int_0^t \int_0^{\infty} \int_{0}^{\infty} \int_0^x \tilde{\omega}(x, y, z) \mathcal{B}(x, y; z, \zeta(s) ) dz dy dx ds.
\end{align*}
Using once more weak convergence $\zeta_n \longrightarrow \zeta ~~~\mbox{in} ~~\mathcal{C}([0,T]; w\text{-}L^1(\mathbb{R}_{>0} ))$, 
and applying \eqref{last}, we conclude
\begin{align*}
\int_0^{\infty} \omega(x) & [ \zeta(t, x)- \zeta^{in}(x) ] dx \\
 = & \lim_{n \to \infty} \int_0^n \omega(x) [ \zeta_n(t, x)- \zeta_n^{in}(x)] dx \nonumber\\
 = & \lim_{n\to \infty} 
  \int_0^t \int_0^{\infty}\int_0^{\infty} \tilde{\omega}(x, y, z)   \mathcal{B}_n(x, y; z, \zeta(s) )  dz dy dx ds  \nonumber\\
=&  \int_0^t \int_0^{\infty} \int_{0}^{\infty} \int_0^x \tilde{\omega}(x, y, z) \mathcal{B}(x, y; z,  \zeta(s) )  dz dy dx ds,
\end{align*}
for every $\omega \in L^{\infty}$. This confirms that $\zeta$ is a weak solution to \eqref{GEDM}--\eqref{Initialdata} in the sense 
of Definition \ref{Definition of Soln} which completes the proof of Theorem \ref{TheoremGEDM2}.
\end{Proof}

\section{Uniqueness}\label{sec:uniq}
In this section, our aim is to prove Theorem~\ref{UniquenessTheorem}. 
\begin{Proof} \emph{of  Theorem~\ref{UniquenessTheorem}}.  We prove it by contradiction. Suppose the solution of \eqref{GEDM}--\eqref{Initialdata} is not unique. 
 Then, there must exist $\zeta$ and $\eta$, two solutions (in a weak sense) to \eqref{GEDM} with $\zeta^{in} = \eta^{in}$. 
Consider $\Delta := \zeta - \eta $ and $L:=\sign(\Delta) $. Let $\omega(x)  := \max\{1, x^{1/2}\}$, then consider the following integral:
\begin{align}\label{Uniqueness1}
& \frac{d}{dt} \int_0^{\infty} \omega(x) |\Delta(t, x)| dx  \nonumber\\
= & \int_0^{\infty} \int_0^{\infty} \int_z^{\infty}  \mathcal{A}(x, y; z) \Theta(t, x, y, z) \{ \zeta(t, x) \zeta(t, y) - \eta(t, x) \eta(t, y)  \}  dx dy dz \nonumber\\
= & \int_0^{\infty} \int_0^{\infty} \int_z^{\infty}  \Theta(t, x, y, z) \mathcal{A}(x, y; z)   \zeta(t, x)   \Delta(t, y)   dx  dy dz \nonumber\\
 & + \int_0^{\infty} \int_0^{\infty} \int_z^{\infty}   \Theta(t, x, y, z)  \mathcal{A}(x, y; z)  \eta(t, y) \Delta(t, x) dx  dy dz,
\end{align}
where 
\begin{align*}
 \Theta(t, x, y, z) := \omega(y + z) L(t, y+z) + \omega(x - z) L(t, x-z) - \omega(x) L(t, x) - \omega(y) L(t, y).
\end{align*}
We can see that 
\begin{align}\label{Bound11}
 \Theta(t, x, y, z) \Delta(t, y) \le \{ \omega(y+z)  + \omega(x-z) + \omega(x) - \omega(y) \} | \Delta(t, y) | ,
\end{align}
and
\begin{align}\label{Bound12}
 \Theta(t, x, y, z) \Delta(t, x) \le \{ \omega(y+z) + \omega(x-z) - \omega(x) + \omega(y) \} | \Delta(t, x) |.
\end{align}
Using \eqref{Bound11} and \eqref{Bound12} we can write \eqref{Uniqueness1} in the form
\begin{align}\label{Uniqueness2}
& \frac{d}{dt} \int_0^{\infty} \omega(x) |\Delta(t, x)| dx  \nonumber\\
\leq   &  \int_0^{\infty} \int_0^{\infty} \int_z^{\infty}   \mathcal{A}(x, y; z) \{ \omega(y+z)  - \omega(y) \} | \Delta(t, y) |    \zeta(t, x)     dx  dy dz \nonumber\\
& + \int_0^{\infty} \int_0^{\infty} \int_z^{\infty}   \mathcal{A}(x, y; z) \{  \omega(x-z) + \omega(x)  \} | \Delta(t, y) |    \zeta(t, x)     dx  dy dz \nonumber\\
& + \int_0^{\infty} \int_0^{\infty} \int_z^{\infty}     \mathcal{A}(x, y; z)  \{ \omega(y+z) +  \omega(y) \} | \Delta(t, x) |  \eta(t, y)  dx  dy dz\nonumber\\
& + \int_0^{\infty} \int_0^{\infty} \int_z^{\infty}     \mathcal{A}(x, y; z) \{  \omega(x-z) - \omega(x)  \} | \Delta(t, x) |  \eta(t, y)  dx  dy dz.
\end{align}
We now consider the following estimates 
\begin{equation*}
\omega(y+z)- \omega(y) =\begin{cases}
0,\           &  \text{if}\ y+z \leq 1, \  y \leq 1, \ z \leq 1, \\
(y+z)^{1/2} -1 \leq z^{1/2} ,\               &  \text{if}\ y+z > 1, \  y \leq 1, \ z \leq 1, \\
(y+z)^{1/2} -1 \leq z^{1/2},\               &  \text{if} \  y \leq  1, \ z > 1, \\
(y+z)^{1/2} -y^{1/2} \leq z^{1/2},\ \       &  \text{if} \  y > 1, \ z > 0, 
\end{cases}   
\end{equation*}
\begin{equation*}
\omega(x-z)+\omega(x) =\begin{cases}
2,\       &  \text{if}\ x \leq 1, z \leq x,\ \\
1+ x^{1/2} \leq 2 x^{1/2},\       &  \text{if}\ x > 1, \  x-z \leq 1, \ z \leq 1,\ \\
1+ x^{1/2} \leq 2 x^{1/2},\       &  \text{if}\ x > 1, \  x-z \leq 1, \  z > 1,\ \\
(x-z)^{1/2}+ x^{1/2} \leq 2x^{1/2},\       &  \text{if}\ x >1, \  x-z > 1,\ 0< z<x,
\end{cases}   
\end{equation*}
\begin{equation*}
\omega(y+z)+ \omega(y) =\begin{cases}
2,\       &  \text{if}\ y+z \leq 1, \  y \leq 1, \ z \leq 1, \\
(y+z)^{1/2} + 1 \leq 3,\       &  \text{if}\ y+z > 1, \  y \leq 1, \ z \leq 1, \\
(y+z)^{1/2} + 1 \leq 2 + z^{1/2},\       &  \text{if}\  y \leq 1, \ z > 1, \\
(y+z)^{1/2} + y^{1/2} \leq 2 y^{1/2} + z^{1/2},\       &  \text{if}\ y > 1,\  z>0,
\end{cases}   
\end{equation*}
and 
\begin{equation*}
\omega(x-z)-\omega(x) = \begin{cases}
0,\       &  \text{if}\ x \leq 1, z \leq x,\ \\
1 - x^{1/2}\leq 0,\       &  \text{if}\ x > 1, \  x-z \leq 1, \ z \leq 1,\ \\
1 - x^{1/2} \leq 0,\       &  \text{if}\ x > 1, \  x-z \leq 1, \  z > 1,\ \\
(x-z)^{1/2} - x^{1/2} \leq 0,\       &  \text{if}\ x >1, \  x-z > 1,\ 0< z<x.
\end{cases}   
\end{equation*}
Using these estimates in \eqref{Uniqueness2} and defining $p(x,y,z):= (1+x)^{1/2} (1+y)^{1/2} \varphi(z) $ we have

\begin{align}\label{Uniqueness3}
& \frac{d}{dt} \int_0^{\infty} \omega(x) |\Delta(t, x)| dx  \nonumber\\
\leq   & A \Bigl(\int_0^{1} \int_{1-z}^{\infty} \int_z^{\infty}+\int_1^{\infty} \int_0^{\infty} \int_z^{\infty}  \Bigr)
p(x,y,z) z^{1/2}  |\Delta(t, y) |    \zeta(t, x)     dx  dy dz \nonumber\\
& + 2A \int_0^{1} \int_0^{\infty} \int_z^{1}  p(x,y,z) | \Delta(t, y) |    \zeta(t, x)     dx  dy dz \nonumber\\
& + 2A\Bigl(\int_0^{\infty} \int_0^{\infty} \int_1^{1+z} + \int_0^{\infty} \int_0^{\infty} \int_{1+z}^{\infty}\Bigr)
p(x,y,z)x^{1/2} | \Delta(t, y) |    \zeta(t, x)     dx  dy dz \nonumber\\
& + 3A \Bigl( \int_0^{1} \int_0^{1-z} \int_z^{\infty} + \int_0^{1} \int_{1-z}^{1} \int_z^{\infty}  \Bigr)
p(x,y,z)| \Delta(t, x) |  \eta(t, y)  dx  dy dz\nonumber\\
& +A \int_1^{\infty} \int_0^{1} \int_z^{\infty}   p(x,y,z)  \{ 2 + z^{1/2} \} | \Delta(t, x) |  \eta(t, y)  dx  dy dz\nonumber\\
& +A \int_0^{\infty} \int_0^{\infty} \int_z^{\infty}   p(x,y,z)\{  2 y^{1/2} + z^{1/2}  \} | \Delta(t, x) |  \eta(t, y)  dx  dy dz.
\end{align}
Remembering that $\varphi\in L^1_{0,1}$ and $\omega(x)  := \max\{1, x^{1/2}\}$ we can estimate \eqref{Uniqueness3} as follows
\begin{align*}\label{Uniqueness4}
 \frac{d}{dt} \int_0^{\infty} \omega(x) |\Delta(t, x)| dx  
\leq   34A\Gamma  \| \varphi \|   \int_{0}^{\infty}    \omega(x)   |\Delta(t, x) |  dx,
\end{align*}
and an application of Gronwall's inequality using $\zeta^{in} = \eta^{in}$ gives
\begin{align*}
 \int_0^{\infty} \omega(x) |\Delta(t, x)| dx \le e^{34 A \Gamma \| \varphi \|  t} \bigg(  \int_0^{\infty} \omega(x) |\Delta(0, x)| dx \bigg) =0.
\end{align*}
Since the function $\omega(x) = \max\{1, x^{1/2}\}$ is positive, this inequality ensures that
\begin{align*}
  \Delta(t, x) =0\ \ \text{almost everywhere},
\end{align*}
i.e $\zeta = \eta $ a.e. This completes the proof of Theorem~\ref{UniquenessTheorem}.
\end{Proof}

\section{Conservation of Mass and Total number of Particles}\label{sec:conserve}
In this section we will establish two results: firstly, we prove that our solution satisfies the conservation of the total number of particles, and then 
we  prove that the total mass is conserved.  The proof of these results follows  ideas similar to those in \cite{Barik:2020, Stewart:1991}.
\begin{Proof} \emph{of Theorem~\ref{TPTheorem}}.  For each $p>0$ let $ \omega (x)  =   \chi_{(0, p]} (x) $ for $x \in (0, \infty )$,
 $\omega_1(x,y):= \omega(x+y)- \omega(x) - \omega(y)$ and $\omega_2 (x, y) := \omega(x-y)- \omega(x) + \omega(y).$  We then obtain
\begin{equation*}
\omega_1 (y, z)  = \begin{cases}
-1,\       &  \text{if}\ y+z \leq p, \  y \leq p, \ z \leq p, \\
-2,\       &  \text{if}\ y+z > p, \  y \leq p, \ z \leq p, \\
-1,\       &  \text{if} \  y \leq p, \ z >p, \\
-1,\       &  \text{if} \  y > p, \ z \leq p, \\
0,\       &  \text{if}\ y > p,\ z > p,
\end{cases}   
\end{equation*}
and 
\begin{equation*}
\omega_2 (x, z) =\begin{cases}
1,\       &  \text{if}\ x \leq p, z \leq x,\ \\
2,\       &  \text{if}\ x > p, \  x-z \leq p, \ z \leq p,\ \\
1,\       &  \text{if}\ x > p, \  x-z \leq p, \  z > p,\ \\
1,\       &  \text{if}\ x >p, \  x-z > p,\ z \leq p,\\
0,\       &  \text{if}\ x>p,\ x-z > p,\ z > p.\\
\end{cases}   
\end{equation*}
Substituting  $\omega_1$ and $\omega_2$ into \eqref{Weak form}, we obtain
\begin{align*}
  \int_0^{p}   [ \zeta(t, x) -   \zeta^{in}(x) ]  
 = &  \int_0^t\Big( - \int_0^p \int_0^{p-z} \int_z^{\infty}  
  - 2   \int_0^p \int_{p-z}^p \int_z^{\infty}  \\
& -   \int_p^{\infty} \int_{0}^{p} \int_{z}^{\infty}   
  - \int_0^{p}  \int_p^{\infty}  \int_z^{\infty}  \\
 & +   \int_0^{p}  \int_0^{\infty}  \int_z^{p}  
 +    2  \int_0^{p}  \int_0^{\infty}  \int_{p}^{p+z}   \\
& +   \int_{p}^{\infty} \! \int_0^{\infty} \! \int_{z}^{p+z}  
 +  \int_0^{p}  \!\int_0^{\infty} \! \int_{p+z}^{\infty}  \Big)  \mathcal{B}(x, y; z, \zeta(s) ) dx  dy  dz  ds.
\end{align*} 
After rearranging the above terms, taking the limit as $p \to \infty$ and then by \eqref{ConditionIntial} we get
\begin{align*}
  \int_0^{\infty}   [ \zeta(t, x)  -   \zeta^{in}(x) ]  dx
 = & \lim_{p \to \infty } \bigg(- \int_0^t  \int_0^{\infty} \int_0^{\infty} \int_z^{\infty}  \ \mathcal{B}(x, y; z, \zeta(s) )   dx dy  dz  ds \nonumber\\
 & -  \int_0^t  \int_0^p \int_{p-z}^p \int_z^{\infty}  \ \mathcal{B}(x, y; z, \zeta(s) )   dx dy  dz  ds \nonumber\\
  &+ \int_0^t  \int_0^{\infty}  \int_0^{\infty}  \int_z^{\infty}   \ \mathcal{B}(x, y; z, \zeta(s) ) dx  dy  dz  ds\nonumber\\
 &+   \int_0^t  \int_0^{p}  \int_0^{\infty}  \int_{p}^{p+z}   \ \mathcal{B}(x, y; z, \zeta(s) ) dx  dy  dz  ds \bigg).
\end{align*} 
Due to \eqref{SumK}, \eqref{Uboundlemma}, \eqref{PhiIntegrable} and \eqref{Convexp1}, it can be seen  
that each integral in the above equation is finite. Consequently, the first and third integrals in the right-hand side
cancel out. Similarly, with the application of \eqref{Convexp3}, \eqref{SumK}, and \eqref{Uboundlemma}, 
the second and fourth integrals tend to $0$ as $p \to \infty $.
This concludes the proof of Theorem~\ref{TPTheorem}.
\end{Proof}

We will now proceed to prove Theorem \ref{MassTheorem}. To achieve this we first prove three auxiliary results.
\begin{lemma}\label{EGDMassLemmaP1}
 Let $ \zeta $ be the weak solution to \eqref{GEDM}--\eqref{Initialdata}. Then, 
\begin{align*}
 \int_0^{p}  x [ \zeta(t, x) -   \zeta^{in}(x) ]  dx 
 =  & - \int_0^t  \int_0^p \int_{p-z}^p \int_z^{\infty} (y+z) \ \mathcal{B}(x, y; z, \zeta(s) )   dx dy  dz  ds \nonumber\\
 &- \int_0^t  \int_p^{\infty} \int_{0}^{p} \int_{z}^{\infty}  y \ \mathcal{B}(x, y; z, \zeta(s) ) dx  dy  dz ds\nonumber\\
 & - \int_0^t \int_0^{p}  \int_p^{\infty}  \int_z^{\infty} z \ \mathcal{B}(x, y; z, \zeta(s) )  dx dy  dz  ds \nonumber\\
 &+ \int_0^t  \int_0^{p}  \int_0^{\infty}  \int_{p}^{p+z}  x \ \mathcal{B}(x, y; z, \zeta(s) ) dx  dy  dz  ds\nonumber\\
& + \int_0^t  \int_{p}^{\infty}  \int_0^{\infty}  \int_{z}^{p+z}  (x-z) \ \mathcal{B}(x, y; z, \zeta(s) ) dx  dy dz  ds\nonumber\\
& + \int_0^t  \int_0^{p}  \int_0^{\infty}  \int_{p+z}^{\infty}  z \ \mathcal{B}(x, y; z, \zeta(s) ) dx  dy  dz  ds.
\end{align*} 
\end{lemma}

\begin{proof}  
Let $ \omega (x)  :=  x \chi_{(0, p]} (x)$, and consider $\omega_1$ and $\omega_2$ as in the proof of
Theorem~\ref{TPTheorem}. Then the proof of this lemma can be readily deduced by substituting the values of $\omega_1$ and $\omega_2$ into \eqref{Weak form}.
\end{proof}

\begin{lemma}\label{EDGMassLemmaP2}
 Let $ \zeta $ be the weak solution to \eqref{GEDM}--\eqref{Initialdata}. Assume $\mathcal{A} $ satisfies \eqref{SumK}. Then, 
\begin{align*}
\int_p^{\infty}   [ \zeta(t, x) -   \zeta^{in}(x) ]  dx 
 = & \int_0^t  \Big(\int_0^{p} \int_{p-z}^{p}  \int_{z}^{\infty} -  \int_p^{\infty} \int_{p}^{\infty}  \int_{z}^{\infty}- \int_0^{p} \int_0^{\infty}   \int_{p}^{p+z} \\
 &\ \qquad   +   \int_p^{\infty} \int_0^{\infty}   \int_{p+z}^{\infty}   \Big) \mathcal{B}(x, y; z, \zeta(s) ) dx  dy  dz  ds,
\end{align*} 
and
\begin{align*}
\lim_{p \to \infty } p \int_p^{\infty}   [ \zeta(t, x) -   \zeta^{in}(x) ]  dx=0.
\end{align*} 
\end{lemma}
\begin{proof} To prove the first equality above set $ \omega (x)  =   \chi_{[p, \infty)} (x) $ for $x \in (0, \infty )$.
Then, for $\omega_1$ and $\omega_2$ as in the proof of of Theorem~\ref{TPTheorem}, we have
\begin{equation*}
\omega_1 (y, z) =\begin{cases}
0,\       &  \text{if}\ y+z < p, \  y < p,\ z < p,\ \\
1,\       &  \text{if}\ y+z \geq  p, \  y < p,\ z<p, \\
0,\       &  \text{if} \  y < p,\ z \geq p, \\
0,\       &  \text{if}\ y \geq p,\ z < p,\\
-1,\       &  \text{if}\ y \geq p,\ z \geq p,
\end{cases}   
\end{equation*}
and 
\begin{equation*}
\omega_2 (x, z) =\begin{cases}
0,\       &  \text{if}\ x < p,\ z < x\\
-1,\       &  \text{if}\ x \geq  p, \  x-z < p,\ z <p, \\
0,\       &  \text{if}\ x \geq  p, \  x-z < p,\ z \geq p, \\
0,\       &  \text{if}\ x \geq  p, \  x-z \geq p,\ z < p, \\
1,\       &  \text{if}\ x-z \geq p,\ \ z \geq p.
\end{cases}   
\end{equation*}
Substituting the values of $\omega_1$ and $\omega_2$ into \eqref{Weak form}, we conclude the result. 
Next, we multiply $p$ and then take the limit as $p \to \infty$ of the first equation in Lemma \ref{EDGMassLemmaP2}  to obtain
\begin{align}\label{Masseqn1}
 &\lim_{p \to \infty} \bigg( p \int_0^t  \int_0^{p} \int_{p-z}^{p}  \int_{z}^{\infty}   \ \mathcal{B}(x, y; z, \zeta(s) )  dx dy  dz  ds\nonumber\\
 &\qquad\quad - \int_0^t  \int_p^{\infty} \int_{p}^{\infty}  \int_{z}^{\infty}   \ \mathcal{B}(x, y; z, \zeta(s) )  dx dy  dz  ds\nonumber\\
&\qquad\quad - \int_0^t  \int_0^{p} \int_0^{\infty}   \int_{p}^{p+z}   \ \mathcal{B}(x, y; z, \zeta(s) ) dx  dy  dz  ds\nonumber\\
&\qquad\quad + \int_0^t  \int_p^{\infty} \int_0^{\infty}   \int_{p+z}^{\infty}   \ \mathcal{B}(x, y; z, \zeta(s) ) dx  dy  dz  ds \bigg) \nonumber\\
= & \lim_{p \to \infty}  \int_p^{\infty}  p  [ \zeta(t, x) - \zeta^{in}(x) ]  dx.
\end{align} 
We can readily deduce from the integrability of the function $x \mapsto x  \zeta^{in}(x)$ and the Lebesgue-dominated convergence theorem that 
\begin{align*}
 & \lim_{p \to \infty}  \bigg|  \int_p^{\infty}  p [ \zeta(t, x) - \zeta^{in}(x) ]  dx \bigg|  \leq  \lim_{p \to \infty}  2 \int_p^{\infty}  x   \zeta^{in}(x)   dx =0,
\end{align*} 
and this completes the proof of the lemma.
\end{proof}

\begin{lemma}\label{EDGMassLemmaP3}
 Let $ \zeta $ be the weak solution to \eqref{GEDM}--\eqref{Initialdata}. Assume $\mathcal{A} $ satisfies \eqref{SumK}. Then, 
\begin{align*}
&(i) \ \   \lim_{p \to \infty } p 
  \int_0^t  \int_p^{\infty} \int_{p}^{\infty}  \int_{z}^{\infty}   \ \mathcal{B}(x, y; z, \zeta(s) )  dx dy  dz  ds  = 0,\\
&(ii)  \ \   \lim_{p \to \infty } p  \int_0^t  \int_p^{\infty} \int_0^{\infty}   \int_{p+z}^{\infty}   \ \mathcal{B}(x, y; z, \zeta(s) ) dx  dy  dz  ds =0.
\end{align*} 
\end{lemma}

\begin{proof} Applying \eqref{Uboundlemma}, \eqref{SumK}, and \eqref{PhiIntegrable}, we have
 \begin{align*}
 \int_p^{\infty} \int_{p}^{\infty}  \int_{z}^{\infty}  p  \ \mathcal{B}(x, y; z, \zeta(s) )  dx dy  dz 
\leq & A \int_p^{\infty} \int_{p}^{\infty}  \int_{p}^{\infty}  p  (x+y) \varphi(z) \zeta(s, x) \zeta(s, y)  dx dy  dz \nonumber\\
\leq & 2A \int_p^{\infty} \int_{p}^{\infty}  \int_{p}^{\infty}   x y \varphi(z) \zeta(s, x) \zeta(s, y)  dx dy  dz \nonumber\\
\leq & 2 A \| \varphi \|_{L^1}   \mathcal{M}_1  \int_{p}^{\infty}   x  \zeta(s, x)   dx. 
 \end{align*}   
 We then apply the Lebesgue-dominated convergence theorem to deduce Lemma \ref{EDGMassLemmaP3} $(i)$.

Similarly, by using \eqref{Uboundlemma}, \eqref{SumK}, and \eqref{PhiIntegrable}, we have
  \begin{align*}
\int_p^{\infty} \int_0^{\infty}   \int_{p+z}^{\infty}   \ p \mathcal{B}(x, y; z, \zeta(s) ) dx  dy  dz 
\leq & A \int_p^{\infty} \int_{p}^{\infty}  \int_{p}^{\infty}  p (x+y) \varphi(z) \zeta(s, x) \zeta(s, y)  dx dy  dz \nonumber\\
\leq & 2 A \int_p^{\infty} \int_{p}^{\infty}  \int_{z}^{\infty}   x y \varphi(z) \zeta(s, x) \zeta(s, y)  dx dy  dz \nonumber\\
\leq & 2 A \| \varphi \|_{L^1}   \mathcal{M}_1  \int_{p}^{\infty}   x  \zeta(s, x)   dx. 
 \end{align*} 
Again applying the Lebesgue-dominated convergence theorem, we deduce Lemma \ref{EDGMassLemmaP3} $(ii)$.
\end{proof}

Applying Lemma \ref{EDGMassLemmaP3} and \eqref{Masseqn1} into Lemma \ref{EDGMassLemmaP2}, we obtain the following result.
\begin{align}\label{EDGP4}
& \lim_{p \to \infty } p   \int_0^t  \int_0^{p} \int_{p-z}^{p}  \int_{z}^{\infty}   \ \mathcal{B}(x, y; z, \zeta(s) )  dx dy  dz  ds  \nonumber\\ 
& \qquad\qquad = \lim_{p \to \infty } p   
  \int_0^t  \int_0^{p} \int_0^{\infty}   \int_{p}^{p+z}   \ \mathcal{B}(x, y; z, \zeta(s) ) dx  dy  dz  ds.
\end{align} 
At this point, we can now complete the proof of Theorem \ref{MassTheorem}.
\begin{Proof} \emph{of Theorem \ref{MassTheorem}}.  On one hand, by applying Lemma \ref{EGDMassLemmaP1}, \eqref{EDGP4}, \eqref{Uboundlemma}, \eqref{SumK}, \eqref{PhiIntegrable} and the Lebesgue-dominated convergence theorem, we can derive
\begin{align*}
  \int_0^{\infty}  x [ \zeta(t, x) -   \zeta^{in}(x) ]  dx 
 =  &\lim_{p \to \infty} \bigg( - \int_0^t  \int_0^p \int_{p-z}^p \int_z^{\infty} (y+z) \ \mathcal{B}(x, y; z, \zeta(s) )   dx dy  dz  ds \nonumber\\
 &\qquad+ \int_0^t  \int_0^{p}  \int_0^{\infty}  \int_{p}^{p+z}  x \ \mathcal{B}(x, y; z, \zeta(s) ) dx  dy  dz  ds \bigg)\nonumber\\
 \leq   & \lim_{p \to \infty} \bigg( - \int_0^t  \int_0^p \int_{p-z}^p \int_z^{\infty} p \ \mathcal{B}(x, y; z, \zeta(s) )   dx dy  dz  ds \nonumber\\
 &\qquad+ \int_0^t  \int_0^{p}  \int_0^{\infty}  \int_{p}^{p+z}  (p+z) \ \mathcal{B}(x, y; z, \zeta(s) ) dx  dy  dz  ds \bigg)\nonumber\\
 = & \lim_{p \to \infty} \int_0^t  \int_0^{p}  \int_0^{\infty}  \int_{p}^{p+z}  z \ \mathcal{B}(x, y; z, \zeta(s) ) dx  dy  dz  ds.
\end{align*} 
By applying \eqref{Uboundlemma}, \eqref{SumK}, \eqref{PhiIntegrable} and \eqref{Convexp0}, it can easily be shown that
\begin{align*}
 \lim_{p \to \infty} \int_0^t  \int_0^{p}  \int_0^{\infty}  \int_{p}^{p+z}  z \ \mathcal{B}(x, y; z, \zeta(s) ) dx  dy  dz  ds=0.
\end{align*} 
Thus, we can deduce 
\begin{align}\label{MassEq1}
  \int_0^{\infty}  x  \zeta(t, x) dx \leq   \int_0^{\infty}  x    \zeta^{in}(x)   dx.
\end{align} 
On the other hand, we have
\begin{align*}
  \int_0^{\infty}  x [ \zeta(t, x) -   \zeta^{in}(x) ]  dx 
 =  &\lim_{p \to \infty} \bigg( - \int_0^t  \int_0^p \int_{p-z}^p \int_z^{\infty} (y+z) \ \mathcal{B}(x, y; z, \zeta(s) )   dx dy  dz  ds \nonumber\\
 &+ \int_0^t  \int_0^{p}  \int_0^{\infty}  \int_{p}^{p+z}  x \ \mathcal{B}(x, y; z, \zeta(s) ) dx  dy  dz  ds \bigg)\nonumber\\
 \geq   & \lim_{p \to \infty} \bigg( - \int_0^t  \int_0^p \int_{p-z}^p \int_z^{\infty} (p+z) \ \mathcal{B}(x, y; z, \zeta(s) )   dx dy  dz  ds \nonumber\\
 &+ \int_0^t  \int_0^{p}  \int_0^{\infty}  \int_{p}^{p+z}  p \ \mathcal{B}(x, y; z, \zeta(s) ) dx  dy  dz  ds \bigg)\nonumber\\
 = & - \lim_{p \to \infty}  \int_0^t  \int_0^p \int_{p-z}^p \int_z^{\infty} z \ \mathcal{B}(x, y; z, \zeta(s) )   dx dy  dz  ds.
\end{align*} 
As in the previous argument, we can readily observe that
 \begin{align*} 
 \lim_{p \to \infty}  \int_0^t  \int_0^p \int_{p-z}^p \int_z^{\infty} z \ \mathcal{B}(x, y; z, \zeta(s) )   dx dy  dz  ds=0.
\end{align*} 
Consequently, we have
\begin{align}\label{MassEq2}
  \int_0^{\infty}  x  \zeta(t, x) dx \geq   \int_0^{\infty}  x    \zeta^{in}(x)   dx.
\end{align} 
Finally, the proof of this theorem can be inferred by using \eqref{MassEq1} and \eqref{MassEq2}.
\end{Proof}

\section*{Acknowledgments}
\noindent
PKB  acknowledges support from grant PID2020-117846GB-I00, the research network RED2022-134784-T, and the \enquote{María de Maeztu} grant CEX2020-001105-M from the Spanish government. Part of this work was done while PKB enjoyed the hospitality of Centre for Mathematical Analysis, Geometry and Dynamical Systems, at  Instituto Superior T\'ecnico, Lisboa.
FPdC, JTP and RS were partially supported by Funda\c{c}\~ao para a Ci\^encia e Tecnologia (Portugal) grants UIDB/04459/2020 and UIDP/04459/2020.



\medskip
Received xxxx 20xx; revised xxxx 20xx.
\medskip
\end{document}